\documentclass[12pt]{article}
\usepackage{amsfonts}
\usepackage{hyperref}

\def\squarebox#1{\hbox to #1{\hfill\vbox to #1{\vfill}}}
\newcommand{\qed}{\hspace*{\fill}
\vbox{\hrule\hbox{\vrule\squarebox{.667em}\vrule}\hrule}\smallskip}
\newtheorem{teorema}{Theorem}[section]
\newtheorem{lema}[teorema]{Lemma}
\newtheorem{corolario}[teorema]{Corollary}
\newtheorem{proposicao}[teorema]{Proposition}
\newtheorem{defi}[teorema]{Definition}
\newenvironment{profe}{\noindent {\bf Proof:}}{\hfill $\qed $ \newline}
\begin{document}

\title{Conditions for Equality between Lyapunov and Morse Decompositions}
\author{Luciana A. Alves\thanks{%
Supported by FAPESP grant n$^{\mathrm{o}}$ 06/60031-3} \thanks{%
Address: Faculdade de Matem\'{a}tica - Universidade Federal de Uberl\^{a}%
ndia. Campus Santa M\^{o}nica, Av. Jo\~{a}o Naves de \'{A}vila, 2121.
38408-100 Uberl\^{a}ndia - MG, Brasil. e-mail: lualves@famat.ufu.br} and
Luiz A. B. San Martin\thanks{%
Supported by CNPq grant no.\ 303755/09-1, FAPESP grant no.\ 2012/18780-0 and
CNPq/Universal grant no 476024/2012-9.} \thanks{%
Address: Imecc - Unicamp, Departamento de Matem\'{a}tica. Rua S\'{e}rgio
Buarque de Holanda, 651, Cidade Universit\'{a}ria Zeferino Vaz. 13083-859
Campinas - SP, Brasil. e-mail: smartin@ime.unicamp.br}}
\maketitle

\begin{abstract}
Let $Q\rightarrow X$ be a continuous principal bundle whose group $G$ is
reductive. A flow $\phi $ of automorphisms of $Q$ endowed with an ergodic
probability measure on the compact base space $X$ induces two decompositions
of the flag bundles associated to $Q$. A continuous one given by the finest
Morse decomposition and a measurable one furnished by the Multiplicative
Ergodic Theorem. The second is contained in the first. In this paper we find
necessary and sufficient conditions so that they coincide. The equality
between the two decompositions implies continuity of the Lyapunov spectra
under pertubations leaving unchanged the flow on the base space.
\end{abstract}

\noindent \textit{AMS 2010 subject classification: Primary: 37H15, \textit{%
37B35, }37B05. Secondary: 22E46, 20M20.}

\noindent \textit{Key words and phrases:} Lyapunov exponents, Multiplicative
Ergodic Theorem, Morse decomposition, Fiber bundles, Semi-simple Lie groups,
Reductive Lie groups, Flag Manifolds.

\section{Introduction}

The purpose of this paper is to give necessary and sufficient conditions for
the equality of Morse and Oseledets decompositions of a continuous flow on a
flag bundle.

We consider a continuous principal bundle $Q\rightarrow X$ with group $G$,
which is assumed to be semi-simple or reductive. A continuous automorphism $%
\phi \in \mathrm{Aut}\left( Q\right) $ of $Q$ defines a discrete-time flow $%
\phi ^{n}$, $n\in \mathbb{Z}$, on $Q$. For instance $Q\rightarrow X$ could
be the bundle of frames of a $d$-dimensional vector bundle $\mathcal{V}%
\rightarrow X$ over $X$, in which case $G$ is the reductive group $\mathrm{Gl%
}\left( d,\mathbb{R}\right) $. Since a linear flow on a vector bundle lifts
to the bundle of frames our set up includes this classical case.

The flow of automorphisms on $Q$ induces a flow on the base space $X$, also
denoted by $\phi $. It also induces flows on bundles having as typical fiber
a space $F$ acted by $G$. Such bundle is built via the associated bundle
construction and is denoted by $Q\times _{G}F$. If there is no risk of
confusion the flows on the associated bundles are denoted by $\phi $ as well.

When $G$ is a reductive group we are specially interested in their flag
manifolds $\mathbb{F}_{\Theta }$, distinguished by the subindex $\Theta $,
which are compact homogeneous spaces of $G$. We write $\mathbb{E}_{\Theta
}=Q\times _{G}\mathbb{F}_{\Theta }$ for the corresponding flag bundle. For
the flow $\phi $ induced on $\mathbb{E}_{\Theta }$, it was proved in \cite%
{smbflow} and \cite{msm} that it has a finest Morse decomposition (under the
mild assumption that the flow on the base space $X$ is chain transitive).
Each Morse component of this finest decomposition meets a fiber of $\mathbb{E%
}_{\Theta }\rightarrow X$ in a algebraic submanifold of
$\mathbb{F}_{\Theta} $. This submanifold is defined as a set of
fixed points for some $g\in G$ acting on $\mathbb{F}_{\Theta }$.
For instance in a projective bundle the fibers of a Morse
component are subspaces, that can be seen as sets of fixed points
on the projective space of diagonalizable matrices (see also
Selgrade \cite{sel} for the Morse decomposition on a projective
bundle). The Morse decomposition is thus described by a continuous
section $\chi _{\mathrm{Mo}}$
of an associated bundle $Q\times _{G}\left( \mathrm{Ad}\left( G\right) H_{%
\mathrm{Mo}}\right) $, whose typical fiber is an adjoint orbit $\mathrm{Ad}%
\left( G\right) H_{\mathrm{Mo}}$ of $G$. Here $H_{\mathrm{Mo}}$ belongs to
the Lie algebra $\mathfrak{g}$ of $G$ and its adjoint $\mathrm{ad}\left( H_{%
\mathrm{Mo}}\right) $ has real eigenvalues. The Morse components are then
built from the section $\chi _{\mathrm{Mo}}$ and the fixed point sets of $%
\exp H_{\mathrm{Mo}}$ on the flag manifolds. (See \cite{msm}, Theorem 7.5.)

On the other hand we also have the Oseledets decomposition, coming
from the Multiplicative Ergodic Theorem (as proved in
\cite{alvsm}). To consider this decomposition it is required a
$\phi $-invariant measure $\nu $ on the base space. If $\nu $ is
ergodic and $\mathrm{supp}\nu =X$ (which provides chain
transitivity on $X$) then the Multiplicative Ergodic Theorem
yields an analogous decomposition to the Morse decomposition that
describes the level
sets of the $\mathfrak{a}$-Lyapunov exponents (see \cite{alvsm} and \cite%
{smsec}). Again there is an adjoint orbit $\mathrm{Ad}\left( G\right) H_{%
\mathrm{Ly}}$ and a section $\chi _{\mathrm{Ly}}$ of the associated bundle $%
Q\times _{G}\left( \mathrm{Ad}\left( G\right) H_{\mathrm{Ly}}\right) $ such
that the Oseledets decomposition is built from $\chi _{\mathrm{Ly}}$ and the
fixed point sets of $\exp H_{\mathrm{Ly}}$. The section $\chi _{\mathrm{Ly}}$
is now only measurable and defined up to a set of $\nu $-measure $0$.

It turns out that any component of the Oseledets decomposition is contained
in a component of the Morse decomposition (see Section \ref{secversusdec}
below). This means that the multiplicities of the eigenvalues of $\mathrm{ad}%
\left( H_{\mathrm{Mo}}\right) $ are larger than those of $\mathrm{ad}\left(
H_{\mathrm{Ly}}\right) $.

In this paper we write down three conditions that together are necessary and
sufficient for both decompositions to coincide (see Section \ref{threecon}).
In this case the Morse decomposition is a continuous extension of the
Oseledets decomposition.

The first of these conditions requires boundedness of the measurable section
$\chi _{\mathrm{Ly}}$, which means that different components of the
Oseledets decomposition do not approach each other. The other two conditions
are about the Oseledets decomposition for the other ergodic measures on $%
\mathrm{supp}\nu =X$. They can be summarized by saying that if $\rho $ is an
ergodic measure then its Oseledets decomposition is finer than the
decomposition for $\nu $.

It is easy to prove that each of the three conditions is necessary. Our main
result is to prove that together they imply equality of the decompositions.

Now we describe the contents of the paper and say some words about other
results that have independent interest.

Sections \ref{secprelim}, \ref{seclyapmorse} and \ref{secperflow} are
preliminary. Section \ref{secprelim} contains notation \ and general facts
about flag manifolds, while in Section \ref{seclyapmorse} we recall the
results of \cite{alvsm}, \cite{smbflow}, \cite{msm} and \cite{smsec} about
Morse decomposition, Morse and Lyapunov spectra and the Multiplicative
Ergodic Theorem on flag bundles. In Section \ref{secperflow} we discuss
briefly flows over periodic orbits that will be needed later.

Section \ref{secmedinv} is devoted to the analysis of ergodic measures on
the flag bundles. We exploit the Krylov-Bogolyubov technique of occupation
measures to see that any Lyapunov exponent coming from the Multiplicative
Ergodic Theorem is an integral over an ergodic measure and conversely.
Combining this with the fact that an ergodic measure charges just one
Oseledets component allow us to introduce what we called attractor and
repeller measures. Later their supports will provide attractor-repeller
pairs on the flag bundles thus relating them to the finest Morse
decomposition.

In Section \ref{secversusdec} we use the attractor and repeller measures to
check that the components of the Oseledets decomposition are indeed
contained in the finest Morse decomposition.

Another tool is developed in Section \ref{lyextgbun}, namely the Lyapunov
exponents of the derivative flow on the tangent space of the fibers of the
flag bundles. The knowledge of these exponents allow to find $\omega $%
-limits in the bundles themselves.

In Section \ref{teclemma} we prove our main technical lemma that furnishes
attractor-repeller pairs on the flag bundles.

In Section \ref{threecon} we state our conditions and prove that they are
necessary. Their sufficience is proved in Section \ref{sufcon}.

In the next two Sections \ref{secunique} and \ref{seciid} we discuss two
cases that go in opposite directions. Namely flows where the base space is
uniquely ergodic (Section \ref{secunique}) and product of i.i.d. sequences.
For a uniquely ergodic base space the second and third conditions are
vacuous and it follows by previous results that the Morse spectrum is a
polyhedron that degenerates to a point if the first condition is satisfied.
In the other side in the i.i.d. case there are plenty of invariant measures
enabling to find examples that violate our second condition. We do that with
the aid of a result by Guivarch'- Raugi \cite{gr}.

Section \ref{seccont} is independent of the rest of the paper. It contains a
result that motivates the study of the equality between Oseledets and Morse
decompositions. We prove that if both decompositions coincide for $\phi $
then the Lyapunov spectrum is continuous under perturbations $\sigma \phi $
of $\phi $ with $\sigma $ varying in gauge group $\mathcal{G}$ of $Q$. This
continuity is a consequence of the differentiability result of \cite{ferrsm}%
. By that result there exists a subset $\Phi _{\mathrm{Mo}}$ of linear
functionals defined from the finest Morse decomposition such that the map $%
\sigma \mapsto \alpha \left( H_{\mathrm{Ly}}\left( \sigma \phi \right)
\right) $ is differentiable with respect to $\sigma $ (at the identity) if $%
\alpha \in \Phi _{\mathrm{Mo}}$, where $H_{\mathrm{Ly}}\left( \sigma \phi
\right) $ is the vector Lyapunov spectrum of $\sigma \phi $. Having equality
of the decompositions we can exploit upper semi-continuity of the spectrum
to prove continuity of $\beta \left( H_{\mathrm{Ly}}\left( \sigma \phi
\right) \right) $ with $\beta $ in a basis that contains $\Phi _{\mathrm{Mo}%
} $.

Finally, we mention that for a linear flow $\phi $ on a vector bundle $%
\mathcal{V}\rightarrow X$ the topological property given by the finest Morse
decomposition of the flow induced on the projective bundle $\mathbb{P}%
\mathcal{V}\rightarrow X$ can be given an analytic characterization via
exponential separation of vector subbundles of $\mathcal{V}$ (see
Colonius-Kliemann \cite{ck}, Chapter 5, and Bonatti-Diaz-Viana \cite{bodivi}%
, Appendix B). In fact, by a theorem by Bronstein-Chernii \cite{brocher}
(quoted from \cite{ck}) the finest Morse decomposition on $\mathbb{P}%
\mathcal{V}$ corresponds to the finest decomposition of $\mathcal{V}$ into
exponentially separated subbundles (see \cite{ck}, Theorem 5.2.10). Hence
our main result gives, in particular, necessary and sufficient conditions
ensuring that the Oseledets decomposition of a vector bundle is
exponentially separated. Also the result of Section \ref{seccont} shows that
if Oseledets decomposition is exponentially separated then the Lyapunov
spectrum changes continuously when $\phi $ is perturbed in such a way that
the flow on the base $X$ is kept fixed.

\section{Flag manifolds}

\label{secprelim}

We explain here our notation about semi-simple (or reductive) Lie groups and
their flag manifolds. We refer to {Knapp} \cite{knp}, {%
Duistermat-Kolk-Varadarajan} \cite{DKV} and {Warner} \cite{Warner}.

Let $\mathfrak{g}$ be a semi-simple non-compact Lie algebra. In order to
make the paper understandable to readers without acquaintance with Lie
Theory we adopted the strategy of defining the notation by writing
explicitely their meanings for the special linear group $\mathrm{Sl}\left( d,%
\mathbb{R}\right) $ and its Lie algebra $\mathfrak{sl}\left( d,\mathbb{R}%
\right) $ (or $\mathrm{Gl}\left( d,\mathbb{R}\right) $ and $\mathfrak{gl}%
\left( d,\mathbb{R}\right) $ in the reductive case). We hope that the reader
with expertise in semi-simple theory will recognize the notation for the
general objects (e.g. $\mathfrak{k}$ is a maximal compact embedded
subalgebra, etc.)

At the Lie algebra level the Cartan decomposition reads $\mathfrak{g}=%
\mathfrak{k}\oplus \mathfrak{s}$, where $\mathfrak{k}=\mathfrak{so}\left(
d\right) $ is the subalgebra of skew-symmetric matrices and $\mathfrak{s}$
is the space of symmetric matrices. The Iwasawa decomposition of the Lie
algebra is $\mathfrak{g}=\mathfrak{k}\oplus \mathfrak{a}\oplus \mathfrak{n}$
where $\mathfrak{a}$ is the subalgebra of diagonal matrices and $\mathfrak{n}
$ is the subalgebra of upper triangular matrices with zeros on the diagonal.

The set of roots is denoted by $\Pi $. These are linear maps $\alpha
_{ij}\in \mathfrak{a}^{*}$, $i\neq j$, defined by $\alpha _{ij}\left(
\mathrm{diag}\{a_{1},\ldots ,a_{d}\}\right) =a_{i}-a_{j}$. The set of
positive roots is $\Pi ^{+}=\{\alpha _{ij}:i<j\}$ and the set of simple
roots is $\Sigma =\{\alpha _{ij}:j=i+1\}$. The root space is $\mathfrak{g}%
_{\alpha }$ ($\mathfrak{g}_{\alpha _{ij}}$ is spanned by the basic matrix $%
E_{ij}$) and
\[
\mathfrak{g}=\mathfrak{m}\oplus \mathfrak{a}\oplus \sum_{\alpha \in \Pi }%
\mathfrak{g}_{\alpha }
\]
where $\mathfrak{m}=\mathfrak{z}_{\mathfrak{k}}\left( \mathfrak{a}\right) =%
\mathfrak{z}\left( \mathfrak{a}\right) \cap \mathfrak{k}$ is the centralizer
of $\mathfrak{a}$ in $\mathfrak{k}$ ($\mathfrak{m}=0$ in $\mathfrak{sl}%
\left( d,\mathbb{R}\right) $). The basic (positive) Weyl chamber is denoted
by
\[
\mathfrak{a}^{+}=\{H\in \mathfrak{a}:\alpha \left( H\right) >0,\alpha \in
\Sigma \}
\]
(cone of diagonal matrices $\mathrm{diag}\{a_{1},\ldots ,a_{d}\}$ satisfying
$a_{1}>\cdots >a_{d}$). Its closure $\mathrm{cl}\mathfrak{a}^{+}$ is formed
by diagonal matrices with decreasing eigenvalues.

At the Lie group level the Cartan decomposition reads $G=KS$, $K=\exp
\mathfrak{k}$ and $S=\exp \mathfrak{s}$ ($K$ is the group $\mathrm{SO}\left(
d\right) $ and $S$ the space of positive definite symmetric matrices in $%
\mathrm{Sl}\left( d,\mathbb{R}\right) $). The Iwasawa decomposition is $%
G=KAN $, $A=\exp \mathfrak{a}$, $N=\exp \mathfrak{n}$. The Cartan
decomposition splits further into the polar decomposition $G=K\left( \mathrm{%
cl}A^{+}\right) K$, $A^{+}=\exp \mathfrak{a}^{+}$.

$M=\mathrm{Cent}_{K}\left( \mathfrak{a}\right) $ is the centralizer of $%
\mathfrak{a}$ in $K$ (diagonal matrices with entries $\pm 1$), $M^{\ast }=%
\mathrm{Norm}_{K}\left( \mathfrak{a}\right) $ is the normalizer of $%
\mathfrak{a}$ in $K$ (signed permutation matrices) and $\mathcal{W}=M^{\ast
}/M$ is the \textit{Weyl group} (for $\mathrm{Sl}\left( d,\mathbb{R}\right) $
it is the group of permutations in $d$ letters, that acts in $\mathfrak{a}$
by permuting the entries of a diagonal matrix).

The (standard) minimal parabolic subalgebra is $\mathfrak{p}=$ $\mathfrak{m}%
\oplus \mathfrak{a}\oplus \mathfrak{n}$ (= upper triangular matrices), and a
general standard parabolic subalgebra $\mathfrak{p}_{\Theta }$ is defined by
a subset $\Theta \subset \Sigma $ as
\[
\mathfrak{p}_{\Theta }=\mathfrak{m}\oplus \mathfrak{a}\oplus \sum_{\alpha
\in \Pi ^{+}}\mathfrak{g}_{\alpha }\oplus \sum_{\alpha \in \langle \Theta
\rangle ^{+}}\mathfrak{g}_{-\alpha },
\]%
where $\langle \Theta \rangle $ is the set of roots spanned (over $\mathbb{Z}
$) by $\Theta $ and $\langle \Theta \rangle ^{+}=\langle \Theta \rangle \cap
\Pi ^{+}$. That is, $\mathfrak{p}_{\Theta }=\mathfrak{p}\oplus \mathfrak{n}%
^{-}(\Theta )$, where $\mathfrak{n}^{\pm }(\Theta )=\sum_{\alpha \in \langle
\Theta \rangle ^{+}}\mathfrak{g}_{\pm \alpha }$.

Alternatively, given $\Theta $, take $H_{\Theta }\in \mathrm{cl}\mathfrak{a}%
^{+}$ such that $\alpha \left( H_{\Theta }\right) =0$, $\alpha \in \Sigma $,
if and only if $\alpha \in \Theta $. Such $H_{\Theta }$ exists and we call
it a characteristic element of $\Theta $. Then $\mathfrak{p}_{\Theta }$ is
the sum of eigenspaces of $\mathrm{ad}\left( H_{\Theta }\right) $ having
eigenvalues $\geq 0$. In $\mathfrak{sl}\left( d,\mathbb{R}\right) $, $%
H_{\Theta }=\mathrm{diag}\left( a_{1},\ldots ,a_{d}\right) $ with $a_{1}\geq
\cdots \geq a_{d}$, where the multiplicities of the eigenvalues is
prescribed by $a_{i}=a_{i+1}$ if $\alpha _{i,i+1}\in \Theta $, that is, $%
\mathfrak{p}_{\Theta }$ is the subalgebra of matrices that are upper
triangular in blocks, whose sizes are the multiplicities of the eigenvalues
of $H_{\Theta }$. When $\Theta $ is empty, $\mathfrak{p}_{\emptyset }$ boils
down to the minimal parabolic subalgebra $\mathfrak{p}$.

Conversely if $H\in \mathrm{cl}\mathfrak{a}^{+}$ then $\Theta _{H}=\{\alpha
\in \Sigma :\alpha \left( H\right) =0\}$ defines a flag manifold $\mathbb{F}%
_{\Theta _{H}}$ (e.g. the Grassmannian $\mathrm{Gr}_{k}\left( d\right) $ is
a flag manifold of $\mathrm{Sl}\left( d,\mathbb{R}\right) $ defined by $H=%
\mathrm{diag}\{a,\ldots ,a,b\ldots ,b\}$ with $\left( n-k\right) a+kb=0$).
For $H_{1},H_{2}\in \mathrm{cl}\mathfrak{a}^{+}$ we say that $H_{1}$ refines
$H_{2}$ in case $\Theta _{H_{1}}\subset \Theta _{H_{2}}$. In $\mathfrak{sl}%
\left( d,\mathbb{R}\right) $ this means that the blocks determined by the
multiplicities of the eigenvalues of $H_{1}$ is contained in the blocks of $%
H_{2}$.

The \textit{parabolic subgroup} $P_{\Theta }$, associated to $\Theta $, is
defined as the normalizer of $\mathfrak{p}_{\Theta }$ in $G$ (as a group of
matrices it has the same block structure as $\mathfrak{p}_{\Theta }$). It
decomposes as $P_{\Theta }=K_{\Theta }AN$, where $K_{\Theta }=\mathrm{Cent}%
_{K}\left( H_{\Theta }\right) $ is the centralizer of $H_{\Theta }$ in $K$.
We usually omit the subscript when $\Theta =\emptyset $ and $P=P_{\emptyset
} $ is the minimal parabolic subgroup.

The \textit{flag manifold} associated to $\Theta $ is the homogeneous space $%
\mathbb{F}_{\Theta }=G/P_{\Theta }$ (just $\mathbb{F}$ when $\Theta
=\emptyset $). If $\Theta _{1}\subset \Theta _{2}$ then the corresponding
parabolic subgroups satisfy $P_{\Theta _{1}}\subset P_{\Theta _{2}}$, so
that there is a canonical fibration $\pi _{\Theta _{2}}^{\Theta _{1}}:%
\mathbb{F}_{\Theta _{1}}\rightarrow \mathbb{F}_{\Theta _{2}}$, given by $%
gP_{\Theta _{1}}\mapsto gP_{\Theta _{2}}$ (just $\pi _{\Theta _{2}}$ if $%
\Theta _{1}=\emptyset $). For the matrix group the flag manifold $\mathbb{F}%
_{\Theta }$ identifies with the manifold of flags of subspaces $V_{1}\subset
\cdots \subset V_{k}$ where the differences $\dim V_{i+1}-\dim V_{i}$ are
the sizes of the blocks defined by $\Theta $ (or rather the diagonal \
matrix $H_{\Theta }$). The projection $\pi _{\Theta _{2}}^{\Theta _{1}}:%
\mathbb{F}_{\Theta _{1}}\rightarrow \mathbb{F}_{\Theta _{2}}$ is defined by
\textquotedblleft forgetting subspaces\textquotedblright .

The concept of \textit{dual flag manifold} is defined as follows: Let $w_{0}$
be the principal involution of $\mathcal{W}$, that is, the only element of $%
\mathcal{W}$ such that $w_{0}\mathfrak{a}^{+}=-\mathfrak{a}^{+}$, and put $%
\iota =-w_{0}$. Then $\iota (\Sigma )=\Sigma $ and for $\Theta \subset
\Sigma $ write $\Theta ^{\ast }=\iota (\Theta )$. Then $\mathbb{F}_{\Theta
^{\ast }}$ is called the flag manifold dual of $\mathbb{F}_{\Theta }$. For
the matrix group the vector subspaces of the flags in $\mathbb{F}_{\Theta
^{\ast }}$ have complementary dimensions to those in $\mathbb{F}_{\Theta }$
(for instance the dual of a Grassmannian $\mathrm{Gr}_{k}\left( d\right) $
is the Grasmmannian $\mathrm{Gr}_{d-k}\left( d\right) $).

We say that two elements $b_{1}\in \mathbb{F}_{\Theta }$ and $b_{2}\in
\mathbb{F}_{\Theta ^{\ast }}$ are transversal if $(b_{1},b_{2})$ belongs to
the unique open $G$-orbit in $\mathbb{F}_{\Theta }\times \mathbb{F}_{\Theta
^{\ast }}$, by the action $g\left( b_{1},b_{2}\right) =\left(
gb_{1},gb_{2}\right) $. For instance $b_{1}\in \mathrm{Gr}_{k}\left(
d\right) $ and $b_{2}\in \mathrm{Gr}_{d-k}\left( d\right) $ are transversal
if and only if they are transversal as subspaces of $\mathbb{R}^{d}$. In
general transversality can be expressed in terms of transversality of
subalgebras of $\mathfrak{g}$ (see e.g. \cite{smmax}). By the very
definition transversality is an open condition and if $g\in G$ then $gb_{1}$
is transversal to $gb_{2}$ if and only if $b_{1}$ is transversal to $b_{2}$.
The following lemma about transversality will be used afterwards.

\begin{lema}
\label{lemseqtransv}Let $b_{n}^{*}$ be a sequence in $\mathbb{F}_{\Theta
^{*}}$ with $\lim b_{n}^{*}=b^{*}$. Suppose that $b\in \mathbb{F}_{\Theta }$
is not transversal to $b^{*}$. Then there exists a sequence $b_{n}\in
\mathbb{F}_{\Theta }$ with $b_{n}$ not transversal to $b_{n}^{*}$ such that $%
\lim b_{n}=b$.
\end{lema}

\begin{profe}
There exists a sequence $k_{n}\in K$ with $b_{n}^{\ast }=k_{n}b^{\ast }$ and
$k_{n}\rightarrow 1$. Since $b$ is not transversal to $b^{\ast }$ it follows
that $k_{n}b$ is not transversal to $b_{n}^{\ast }$. Hence, $b_{n}=k_{n}b$
is the required sequence.
\end{profe}

We consider now the fixed point set of the action of $h=\exp H$, $H\in
\mathrm{cl}\mathfrak{a}^{+}$, on a flag manifold $\mathbb{F}_{\Theta }$.
Look first at the example of the projective space $\mathbb{R}P^{d-1}$. The
fixed point set is the union of the eigenspaces of $h$. The eigenspace
associated to the biggest eigenvalue is the only attractor (for the
iterations $h^{n}$) that has an open and dense stable manifold. The same way
the eigenspace of the smallest eigenvalue is the unique repeller with open
and dense unstable manifold.

In general the flow defined by $\exp tH$ is gradient in any flag manifold $%
\mathbb{F}_{\Theta }$ (see \cite{DKV}). Its fixed point set is given the
union of the orbits
\[
Z_{H}\cdot wb_{\Theta }=K_{H}\cdot wb_{\Theta }\qquad \,w\in \mathcal{W},
\]%
where $b_{\Theta }=P_{\Theta }$ the origin of $\mathbb{F}_{\Theta
}=G/P_{\Theta }$, $Z_{H}=\left\{ g\in G:\,\mathrm{Ad}(g)H=H\right\} $, $%
K_{H}=Z_{H}\cap K$ and $w$ runs through the Weyl group $\mathcal{W}$. We
write $\mathrm{fix}_{\Theta }(H,w)=Z_{H}\cdot wb_{\Theta }$ and refer to it
as the set of $H$-fixed points of type\textit{\ }$w$. In addition, $\mathrm{%
fix}_{\Theta }(H,1)$ is the only attractor while $\mathrm{fix}_{\Theta
}(H,w_{0})$ is the unique repeller, where $w_{0}$ the principal involution
of $\mathcal{W}$.

For the stable and unstable sets of $\mathrm{fix}_{\Theta }(H,w)$ let $%
\Theta _{H}=\{\alpha \in \Sigma :\alpha \left( H\right) =0\}$ and consider
the nilpotent subalgebras
\[
\mathfrak{n}_{H}^{+}=\sum_{\alpha \in \Pi ^{+}\backslash \langle \Theta
_{H}\rangle }\mathfrak{g}_{\alpha }\qquad \mathfrak{n}_{H}^{-}=\sum_{\alpha
\in \Pi ^{+}\backslash \langle \Theta _{H}\rangle }\mathfrak{g}_{-\alpha }
\]%
and the connected subgroups $N_{H}^{\pm }=\exp \mathfrak{n}_{H}^{\pm }$. Put%
\[
\mathrm{st}_{\Theta }(H,w)=N_{H}^{-}K_{H}\cdot wb_{\Theta }\qquad \mathrm{un}%
_{\Theta }(H,w)=N_{H}^{+}K_{H}\cdot wb_{\Theta }.
\]%
Then $\mathrm{st}_{\Theta }(H,w)$ and $\mathrm{un}_{\Theta }(H,w)$ are the
stable and unstable sets of $\mathrm{fix}_{\Theta }(H,w)$, \ respectively.

More generally if $D=\mathrm{Ad}(g)H$, $g\in G$ and $H\in \mathrm{cl}%
\mathfrak{a}^{+}$, then the dynamics of $\exp tD$ is conjugate under $g$ to
the dynamics of $\exp tH$. Hence $\mathrm{fix}_{\Theta }(D,w)=g\cdot \mathrm{%
fix}_{\Theta }(H,w)$, $\mathrm{st}_{\Theta }(D,w)=g\cdot \mathrm{st}_{\Theta
}(H,w)$ and $\mathrm{un}_{\Theta }(D,w)=g\cdot \mathrm{un}_{\Theta }(H,w)$.
It follows that
\[
\mathrm{st}_{\Theta }(D,w)=P_{D}^{-}\cdot gwb_{\Theta }\qquad \mathrm{un}%
_{\Theta }(H,w)=P_{D}^{+}\cdot gwb_{\Theta },
\]%
where $P_{D}^{\pm }=g{N}_{H}^{\pm }K_{H}g^{-1}=N_{D}^{\pm }K_{D}$, and $%
N_{D}^{\pm }=gN_{H}^{\pm }g^{-1}$ and $K_{D}=gK_{H}g^{-1}$.

If $H_{1}$ refines $H_{2}$ then the centralizers satisfy $Z_{H_{1}}\subset
Z_{H_{2}}$ hence the fixed point set of $\exp H_{1}$ in a flag manifold $%
\mathbb{F}_{\Theta }$ is contained in the fixed point set of $\exp H_{2}$.

The following lemma shows that we can control the inclusion of fixed point
sets for different elements by looking at the attractor and repeller fixed
point sets in the right flag manifolds.

\begin{lema}
\label{lemfixinclu}Suppose $H_{1}$ refines $H_{2}$, take $S\in \mathrm{Ad}%
\left( G\right) H_{1}$ and $T\in \mathrm{Ad}\left( G\right) H_{2}$, and put $%
s=\exp S$, $t=\exp T$. Suppose that $\mathrm{att}_{\Theta \left(
H_{1}\right) }\left( s\right) \subset \mathrm{att}_{\Theta \left(
H_{1}\right) }\left( t\right) $ and $\mathrm{rp}_{\Theta \left( H_{1}\right)
^{*}}\left( s\right) \subset \mathrm{rp}_{\Theta \left( H_{1}\right)
^{*}}\left( t\right) $. Then the fixed point set of $s$ in any flag manifold
is contained in the fixed point set of $t$.\newline
Moreover the fixed points are the same in case these attractor and repeller
fixed points coincide.
\end{lema}

\begin{profe}
If we identify $\mathrm{Ad}\left( G\right) H_{1}$ with the open orbit in $%
\mathbb{F}_{\Theta \left( H_{1}\right) }\times \mathbb{F}_{\Theta \left(
H_{1}\right) ^{\ast }}$ then $S$ is identified to the pair $\left( \mathrm{%
att}_{\Theta \left( H_{1}\right) }\left( s\right) ,\mathrm{rp}_{\Theta
\left( H_{1}\right) }\left( s\right) \right) $. The same way $T$ is
identified to the pair $\left( \mathrm{att}_{\Theta \left( H_{2}\right)
}\left( t\right) ,\mathrm{rp}_{\Theta \left( H_{2}\right) }\left( t\right)
\right) \in \mathbb{F}_{\Theta \left( H_{2}\right) }\times \mathbb{F}%
_{\Theta \left( H_{2}\right) ^{\ast }}$. Now, since $H_{1}$ refines $H_{2}$
there are  fibrations $p:\mathrm{Ad}\left( G\right) H_{1}\rightarrow \mathrm{%
Ad}\left( G\right) H_{2}$, $\pi _{1}:\mathbb{F}_{\Theta \left( H_{1}\right)
}\rightarrow \mathbb{F}_{\Theta \left( H_{2}\right) }$ and $\pi _{2}:\mathbb{%
F}_{\Theta \left( H_{1}\right) ^{\ast }}\rightarrow \mathbb{F}_{\Theta
\left( H_{2}\right) ^{\ast }}$ with the equalities $\pi _{1}\left( \mathrm{%
att}_{\Theta \left( H_{1}\right) }\left( t\right) \right) =\mathrm{att}%
_{\Theta \left( H_{2}\right) }\left( t\right) $ and $\pi _{2}\left( \mathrm{%
att}_{\Theta \left( H_{1}\right) }\left( t\right) \right) =\mathrm{att}%
_{\Theta \left( H_{2}\right) }\left( t\right) $, we have $p\left( S\right) =T
$. This means there exists $g\in G$ such that $S=\mathrm{Ad}\left( g\right)
H_{1}$ and $T=\mathrm{Ad}\left( g\right) H_{2}$. Hence the fixed point set
of $s$ (respectively $t$) in a flag manifold $\mathbb{F}_{\Theta }$ is the
image under $g$ of the fixed point set of $\exp H_{1}$ (respectively $\exp
H_{2}$), which implies the lemma.
\end{profe}

\section{Lyapunov and Morse spectra and decompositions\label{seclyapmorse}}

From now on we consider a discrete-time continuous flow $\phi _{n}$ on a
continuous principal bundle $\left( Q,X,G\right) $, where the base space $X$
is a compact metric space endowed with an ergodic invariant measure $\nu $
with $\mathrm{supp}\nu =X$. The structural group $G$ is assumed to be
semi-simple and noncompact, or slightly more generally, $G$ is reductive
with noncompact semi-simple component. We fix once and for all a maximal
compact subgroup $K\subset G$ and a $K$-subbundle $R\subset Q$. (For a
bundle of frames of a vector bundle $\mathcal{V}\rightarrow X$ this amounts
to the choice of a Riemannian metric on $\mathcal{V}$. In case of a trivial
bundle $Q=X\times G$ the reduction is $R=X\times K$.)

The Iwasawa ($G=KAN$) and Cartan ($G=KS$) decompositions of $G$ yield
decompositions of $Q=R\times AN$ and $Q=R\times S$ by writing $q\in Q$ as
\[
q=r\cdot hn\qquad \mathrm{and}\qquad q=r\cdot s
\]%
$r\in R$, $hn\in AN$ and $s\in S$. In what follows we write for $q\in Q$,
\[
\mathsf{a}\left( q\right) =\log \mathsf{A}\left( q\right) \in \mathfrak{a}
\]%
where $\mathsf{A}\left( q\right) $ is the projection onto $A$ against the
Iwasawa decomposition. Also we write $\mathsf{S}:Q\rightarrow S$ as the
projection onto $S$ of $Q=R\times S$. By the polar decomposition $G=K(%
\mathrm{cl}A^{+})K$ we get a map $\mathsf{A}^{+}:Q\rightarrow \mathrm{cl}%
A^{+}$ by $\mathsf{S}\left( q\right) =k\mathsf{A}^{+}\left( q\right) k^{-1}$%
, $k\in K$. We write
\[
\mathsf{a}^{+}(q)=\log \mathsf{A}^{+}(q)\in \mathrm{cl}{\mathfrak{a}}^{+}.
\]

Now the flow $\phi _{n}$ on $Q$ induces a flow $\phi _{n}^{R}$ on $R$ by
declaring $\phi _{n}^{R}\left( r\right) $ to be the projection of $\phi _{n}
$ onto $R$ against the decomposition \ $Q=R\times AN$ ($\phi _{n}^{R}$ is
indeed a flow because $AN$ is a group.) The projections $\mathsf{a}$ and $%
\mathsf{a}^{+}$ define maps (denoted by the same letters) $\mathsf{a}:%
\mathbb{Z}\times R\rightarrow \mathfrak{a}$ and $\mathsf{a}^{+}:\mathbb{Z}%
\times R\rightarrow \mathrm{cl}\mathfrak{a}^{+}$ by
\[
\mathsf{a}\left( n,r\right) =\mathsf{a}\left( \phi _{n}\left( r\right)
\right) \qquad \mathrm{and}\qquad \mathsf{a}^{+}\left( n,r\right) =\mathsf{a}%
^{+}\left( \phi _{n}\left( r\right) \right) .
\]

It turns out that $\mathsf{a}\left( n,r\right) $ is an additive cocycle over
$\phi _{n}^{R}$. This cocycle factors to a cocycle (also denoted by $\mathsf{%
a}$) over the flow induced on $\mathbb{E}=Q\times _{G}\mathbb{F} $, a associated bundle of $Q$ with typical fiber %
the maximal flag manifold $\mathbb{F}$. The $\mathfrak{a}$-Lyapunov
exponent of $\phi _{n}$ in the direction of $\xi \in \mathbb{E}$ is defined
by
\[
\lambda \left( \xi \right) =\lim_{k\rightarrow +\infty }\frac{1}{k}\mathsf{a}%
\left( k,\xi \right) \in \mathfrak{a}\qquad \xi \in \mathbb{E}.
\]%
The polar exponent is defined by
\[
H_{\phi }\left( r\right) =\lim_{k\rightarrow +\infty }\frac{1}{k}\mathsf{a}%
^{+}\left( k,r\right) \in \mathrm{cl}\mathfrak{a}^{+}\qquad r\in R.
\]%
It turns out that $H_{\phi }\left( r\right) $ is constant along the fibers
of $R$ (when it exists) so is written $H_{\phi }\left( x\right) $, $x\in X$.
The existence of these limits is ensured by the

\vspace{12pt}%

\noindent%
\textbf{Multiplicative Ergodic Theorem} (\cite{alvsm}): The polar exponent $%
H_{\phi }\left( x\right) $ exists for $x$ in a set of total measure $\Omega $%
. Assume that $\nu $ is ergodic. Then $H_{\phi }\left( \cdot \right) $ is
almost surely equals to the constant $H_{\mathrm{Ly}}=H_{\mathrm{Ly}}\left(
\nu \right) \in \mathrm{cl}\mathfrak{a}^{+}$. Put $\mathbb{E}_{\Omega }=\pi
^{-1}\left( \Omega \right) $ where $\pi :\mathbb{E}\rightarrow X$ is the
projection. Then,

\begin{enumerate}
\item $\lambda \left( \xi \right) $ exists for every $\xi \in \mathbb{E}%
_{\Omega }$ and the map $\lambda :\mathbb{E}_{\Omega }\rightarrow \mathfrak{a%
}$ assume values in the finite set $\{wH_{\mathrm{Ly}}:w\in \mathcal{W}\}$.

\item There exists a measurable section $\chi _{\mathrm{Ly}}$ of the bundle $%
Q\times _{G}\mathrm{Ad}\left( G\right) \left( H_{\mathrm{Ly}}\right) $,
defined on $\Omega $, such that $\lambda \left( \xi \right) =w^{-1}H_{%
\mathrm{Ly}}$ if $\xi \in \mathrm{st}(\chi _{\mathrm{Ly}}\left( x\right) ,w)$%
, $x=\pi \left( \xi \right) $.
\end{enumerate}

(To be rigorous the stable set \ $\mathrm{st}(\chi _{\mathrm{Ly}}\left(
x\right) ,w)$, simply denoted by $\mathrm{st}(x,w)$, must be defined using the formalism of fiber bundles. If $%
Q=X\times G$ is trivial then $\chi _{\mathrm{Ly}}:X\rightarrow \mathrm{Ad}%
\left( G\right) \left( H_{\mathrm{Ly}}\right) $ and $\mathrm{st}(\chi _{%
\mathrm{Ly}}\left( x\right) ,w)$ is the stable set discussed in the last
section.)

We write $\mathrm{st}(w)$ for the union of the sets $\mathrm{st}(\chi _{%
\mathrm{Ly}}\left( x\right) ,w)$ with $x$ running through $\Omega $. The
same way we let $\mathrm{fix}(w)$ be the union of the fixed point sets $%
\mathrm{fix}(\chi _{\mathrm{Ly}}\left( x\right) ,w)$.

By analogy with the multiplicative ergodic theorem on vector bundles the
union of the sets $\mathrm{fix}(w)$, $w\in \mathcal{W}$, is called the
Oseledets decomposition of $\mathbb{E}$. These sets project to a partial
flag bundle $\mathbb{E}_{\Theta }$ to fixed point sets $\mathrm{fix}_{\Theta
}(w)$ that form the Oseledets decomposition of $\mathbb{E}_{\Theta }$.

To the exponent $H_{\mathrm{Ly}}\left( \nu \right) \in \mathrm{cl}\mathfrak{a%
}^{+}$ we associate the subset of the simple system of roots
\[
\Theta _{\mathrm{Ly}}=\Theta _{\mathrm{Ly}}\left( \nu \right) =\{\alpha \in
\Sigma :\alpha \left( H_{\mathrm{Ly}}\left( \nu \right) \right) =0\}.
\]%
The corresponding flag manifold $\mathbb{F}_{\Theta _{\mathrm{Ly}}}$ and
flag bundle $\mathbb{E}_{\Theta _{\mathrm{Ly}}}=Q\times _{G}\mathbb{F}%
_{\Theta _{\mathrm{Ly}}}$ play a proeminent role in the proofs. (For a
linear flow on a vector bundle $\mathbb{F}_{\Theta _{\mathrm{Ly}}}$ is the
manifold of flags $\left( V_{1}\subset \cdots \subset V_{k}\right) $ of
subspaces of $\mathbb{R}^{d}$ having the same dimensions as the subspaces of
the Oseledets splitting when the Lyapunov spectrum is ordered decreasingly.)
We refer to $\mathbb{F}_{\Theta _{\mathrm{Ly}}}$ as the flag type of $\phi $
with respect to $\nu $.

As another remark we mention that the section $\chi _{\mathrm{Ly}}$ yields
(actually is built from) two sections $\xi $ and $\xi ^{\ast }$ of the flag
bundles $\mathbb{F}_{\Theta _{\mathrm{Ly}} }$ and $\mathbb{%
F}_{\Theta _{\mathrm{Ly}}^{\ast }}$, respectively. Their images are defined
from level sets of Lyapunov exponents and hence are measurable (see \cite%
{alvsm}, Section 7.1).

On the other hand there are continuous decompositions of the flag bundles
(defined the same way as sets of fixed points) obtained by working out the
concept of Morse decomposition of the flows on the bundles (see Conley \cite%
{con} and Colonius-Kliemann \cite{ck}). It was proved in \cite{smbflow} and
\cite{msm} that if the flow on the base space is chain transitive then the
flow on any flag bundle $\mathbb{E}_{\Theta }$ admits a finest Morse
decomposition with Morse sets $\mathcal{M}\left( w\right) $, also
parametrized by $w\in \mathcal{W}$. Analogous to the Oseledets decomposition
the Morse sets are built as fixed point sets defined by a continuous section
of an adjoint bundle $\chi _{\mathrm{Mo}}:X\rightarrow Q\times _{G}\mathrm{Ad%
}\left( G\right) H_{\mathrm{Mo}}$, where $H_{\mathrm{Mo}}\in \mathrm{cl}%
\mathfrak{a}^{+}$ as well. There is just one attractor Morse component which
is given by $\mathcal{M}^{+}=\mathcal{M}\left( 1\right) $. There is a unique repeller
component as well which is $\mathcal{M}^{-}=\mathcal{M}\left( w_{0}\right) $, where $w_{0}$
is the principal involution.

The assumption that the invariant measure $\nu $ is ergodic with support $%
\mathrm{supp}\nu =X$ implies chain transitivity on $X$.

We write
\[
\Theta _{\mathrm{Mo}}=\Theta _{\mathrm{Mo}}\left( \phi \right) =\{\alpha \in
\Sigma :\alpha \left( H_{\mathrm{Mo}}\right) =0\}
\]%
and refer to $\mathbb{F}_{\Theta _{\mathrm{Mo}}}$ as the flag type of $\phi $
(with respect to the Morse decomposition).

The spectral counterpart of the Morse decomposition is the Morse spectrum
associated to the cocycle $\mathsf{a}\left( n,\xi \right) $. This spectrum
was originally defined by Colonius-Kliemann \cite{ck} for a flow on a vector
bundle and extended to flag bundles (and vector valued cocycles) in \cite%
{smsec}. By the results of \cite{smsec}, each Morse set $\mathcal{M}\left(
w\right) $ has a Morse spectrum $\Lambda _{\mathrm{Mo}}\left( w\right) $
which is a compact convex subset of $\mathfrak{a}$ and contains any $%
\mathfrak{a}$-Lyapunov exponent $\lambda \left( \xi \right) $, $\xi \in
\mathcal{M}\left( w\right) $. The attractor Morse component is given by the
identity $1\in \mathcal{W}$ and we write $\Lambda _{\mathrm{Mo}}=\Lambda _{%
\mathrm{Mo}}\left( 1\right) $, which is the only Morse spectrum meeting $%
\mathrm{cl}\mathfrak{a}^{+}$. The Morse spectrum $\Lambda _{\mathrm{Mo}}$
satisfies the following properties:

\begin{enumerate}
\item $\Lambda _{\mathrm{Mo}}$ is invariant under the group $\mathcal{W}%
_{\Theta _{\mathrm{Mo}}}$ generated by reflections with respect to the roots
$\alpha \in \Theta _{\mathrm{Mo}}$. (See \cite{smsec}, Theorem 8.3.)

\item $\alpha \left( H\right) >0$ if $H\in \Lambda _{\mathrm{Mo}}$ and $%
\alpha $ is a positive root that does not belong to the set $\langle \Theta
\rangle ^{+}$ spanned by $\Theta $. (See \cite{smsec}, Corollary 7.4.)
\end{enumerate}

By the last statement $\alpha \left( H_{\mathrm{Ly}}\right) >0$ if $\alpha $
is a simple root outside $\Theta _{\mathrm{Mo}}$ because $H_{\mathrm{Ly}}\in
\Lambda _{\mathrm{Mo}}$. Hence $\alpha \notin \Theta _{\mathrm{Ly}}$ by
definition of $\Theta _{\mathrm{Ly}}$. It follows that $\Theta _{\mathrm{Ly}%
}\subset \Theta _{\mathrm{Mo}}$. Below in Section \ref{secversusdec} we
improve this statement by proving, with the aid of invariant measures on
flag bundles, that the Oseledets decomposition is contained in the Morse
decomposition.

Our objective is to find necessary and sufficient conditions ensuring that $%
\Theta _{\mathrm{Ly}}=\Theta _{\mathrm{Mo}}$, and hence that the Oseledets
decomposition coincides with the Morse decomposition.

\section{Flows over periodic orbits}

\label{secperflow}

Before proceding let us recall the case where the base space is a single
periodic orbit $X=\{x_{0},\ldots ,x_{\omega -1}\}$ of period $\omega $, that
will be used later to reduce some arguments to nonperiodic orbits.

In the periodic case we have $\Theta _{\mathrm{Ly}}=\Theta _{\mathrm{Mo}}$
since, as is well known, the asymptotics depend ultimately on iterations of
a fixed element in the group $G$. Here the principal bundle is $Q=X\times G$
and the flow is given by
\[
\phi \left( x_{i},h\right) =\left( x_{i+1\left( \mathrm{mod}\omega \right)
},A\left( x_{i}\right) h\right)
\]%
for a map $A:X\rightarrow G$, so that
\[
\phi ^{n}\left( x_{i},h\right) =\left( x_{i+n\left( \mathrm{mod}\omega
\right) },g_{n,i}h\right)
\]%
where $g_{n,i}=A\left( x_{i+n-1\left( \mathrm{mod}\omega \right) }\right)
\cdots A\left( x_{i+1}\right) A\left( x_{i}\right) $. We have $%
g_{n+m,i}=g_{n,i+m\left( \mathrm{mod}\omega \right) }g_{m,i}$ so that $%
g_{k\omega ,i}=g_{\omega ,i}^{k}$. Hence the asymptotics of an orbit
starting at a point above $x_{i}$ is dictated by the iterations of the
action of $g_{\omega ,i}$. The iterations for the action of a fixed $g\in G$
on the flag manifolds, as well as the continuous time version of
periodicity, were studied by Ferraiol-Patr\~{a}o-Seco \cite{fepase}. Let $%
g_{n,i}=u_{n,i}h_{n,i}x_{n,i}$ be the Jordan decomposition of $g_{n,i}$ with
$u_{n,i}$, $h_{n,i}$ and $x_{n,i}$ elyptic, hyperbolic and unipotent
respectively. There is a choice of an Iwasawa decomposition $G=KAN$ such
that $u_{n,i}\in K$, $h_{n,i}\in A$ and $x_{n,i}\in N$. It follows that the
Lyapunov spectrum is given by $\log h_{\omega ,i}$, which is the same for
any $i=0,\ldots ,\omega -1$ (because $g_{\omega ,i+1}=A\left( x_{i}\right)
g_{\omega ,i}A\left( x_{i}\right) ^{-1}$). Also, as proved in \cite{fepase}
the Morse decomposition is given by the fixed point sets of $h_{\omega ,i}$.
Hence $\Theta _{\mathrm{Ly}}=\Theta _{\mathrm{Mo}}$.

\section{Invariant measures on the bundles and $\mathfrak{a}$-Lyapunov
exponents}

\label{secmedinv}

Let $\mu $ be an invariant measure for the flow on the maximal flag bundle $%
\pi :\mathbb{E}\rightarrow X$. Then the integral
\[
\int qd\mu \qquad \,q(\xi )=\mathsf{a}(1,\xi )
\]%
is an $\mathfrak{a}$-Lyapunov exponent for the cocycle $\mathsf{a}\left(
n,\xi \right) $ (see \cite{smsec}). On the other hand by applying the
Multiplicative Ergodic Theorem to an invariant measure $\nu $ on the base
space we obtain $\mathfrak{a}$-Lyapunov exponents, which we call regular
Lyapunov exponents with respect to $\nu $ (because they are obtained as
limits of sequences in $\mathfrak{a}$ which in turn comes  from regular
sequences in $G$, see \cite{alvsm}).

In this section we show that these Lyapunov exponents coincide. Namely, if $%
\nu $ is ergodic measure on $X$ then any of its $\mathfrak{a}$-Lyapunov
exponents is an integral over an ergodic measure $\mu $ that projects onto $%
\nu $, i.e., $\pi _{\ast }\mu =\nu $, and conversely any such integral is a
regular Lyapunov exponent.

Fix an ergodic invariant measure $\nu $ on the base space and let $\Omega
\subset X$ be the set of $\nu $-total measure given by the Multiplicative
Ergodic Theorem (as proved in \cite{alvsm}). Recall that
\[
\pi ^{-1}\left( \Omega \right) =\dot{\bigcup_{w\in \mathcal{W}_{\Theta _{%
\mathrm{Ly}}}\backslash \mathcal{W}}}\mathrm{st}(w)
\]%
and $\lambda \left( \xi \right) =w^{-1}H_{\mathrm{Ly}}$ if $\xi \in \mathrm{%
st}(w)$, where $H_{\mathrm{Ly}}$ is the polar exponent with respect to $\nu $%
.

\begin{proposicao}
\label{propmedinvconjest}Let $\mu $ be an ergodic measure on $\mathbb{E}$
that projects onto $\nu $. Then there exists $w\in \mathcal{W}$ such that $%
\mu (\mathrm{st}(w))=1$ and $\mu \left( \mathrm{st}(w^{\prime })\right) =0$
if $\mathcal{W}_{\Theta _{\mathrm{Ly}}}w\neq \mathcal{W}_{\Theta _{\mathrm{Ly%
}}}w^{\prime }$. In this case
\[
\int qd\mu =w^{-1}H_{\mathrm{Ly}}.
\]
\end{proposicao}

\begin{profe}
By ergodicity of $\mu $ and the ergodic theorem applied to $\mu $ and $%
q\left( \xi \right) =\mathsf{a}(1,\xi )$, there exists a measurable set $%
\mathcal{I}\subset \mathbb{E}$ with $\mu (\mathcal{I})=1$ and
\[
\lambda (\xi )=\lim_{k\rightarrow \infty }\frac{1}{k}\mathsf{a}(k,\xi )=\int
qd\mu \qquad \xi \in \mathcal{I}.
\]%
Now $\mu (\pi ^{-1}(\Omega )\cap \mathcal{I})=1$ and $\pi ^{-1}(\Omega )\cap
\mathcal{I}$ is the disjoint union of the sets $\mathrm{st}(w)\cap \mathcal{I%
}$. In each $\mathrm{st}(w)\cap \mathcal{I}$, $w\in \mathcal{W}$, $\lambda $
is defined and is a constant equal to $w^{-1}H_{\mathrm{Ly}}$. Since $%
\lambda $ is constant on $\mathcal{I}$, it follows that $\pi ^{-1}(\Omega
)\cap \mathcal{I}\subset \mathrm{st}\left( w\right) $, for some $w\in
\mathcal{W}$. Then for any $\xi \in \mathcal{I}$,
\[
\int qd\mu =\lambda (\xi )=w^{-1}H_{\mathrm{Ly}}.
\]%
Finally, $\mu \left( \mathrm{st}(w)\right) \geq \mu \left( \pi ^{-1}(\Omega
)\cap \mathcal{I}\right) =1$, which implies that $\mu \left( \mathrm{st}%
(w^{\prime })\right) =0$ if $\mathrm{st}(w^{\prime })\neq \mathrm{st}(w)$
that is if $W_{\Theta _{\mathrm{Ly}}}w\neq W_{\Theta _{\mathrm{Ly}%
}}w^{\prime }$.
\end{profe}

\begin{corolario}
Let $\Lambda _{\mathrm{Mo}}\left( w\right) \subset \mathfrak{a}$ be the
Morse spectrum of the Morse set $\mathcal{M}\left( w\right) $. Then the
extremal points of the compact convex set $\Lambda _{\mathrm{Mo}}\left(
w\right) $ are regular Lyapunov exponents for ergodic measures on the base
space.
\end{corolario}

\begin{profe}
In fact it was proved \cite{smsec} (see Theorem 3.2(6)) that any extremal
point of $\Lambda _{\mathrm{Mo}}\left( w\right) $ is an integral $\int qd\mu
$ with respect to an ergodic measure $\mu $ on $\mathbb{E}$. (See also \cite%
{ck}, Lemma 5.4.10.)
\end{profe}

The converse to the above proposition says that any regular Lyapunov
exponent is the integral of $q$ with respect to some ergodic measure
projecting onto $\nu $. In order to prove the converse we recall the
Krylov-Bogolyubov procedure of constructing invariant measures as occupation
measures (see e.g. \cite{ck}). Let $\psi _{n}$, $n\in \mathbb{Z}$, be a flow
on a compact metric space $Y$. Then this means that
\[
\left( L_{n,x}f\right) \left( x\right) =\frac{1}{n}\sum_{k=0}^{n-1}f\left(
\psi _{k}x\right) ,\qquad x\in Y,
\]%
define linear maps on the space $C_{0}\left( Y\right) $ of continuous
functions, and hence Borel probability measures $\rho _{n}$. An accumulation
point $\rho _{x}=\lim_{k}\rho _{n_{k}}$ is called an (invariant) occupation
measure. When the limit $\widetilde{f}\left( x\right) =\lim_{n}\frac{1}{n}%
\sum_{k=0}^{n-1}f\left( \phi _{k}x\right) $ exists it is an integral $%
\widetilde{f}\left( x\right) =\int f\left( y\right) \mu _{x}\left( dy\right)
$ with respect to an occupation measure. The following properties will be
used below:

\begin{enumerate}
\item Let $\rho $ be an ergodic probability measure on $Y$. Then for $\rho $%
-almost every $y\in Y$, any occupation measure $\rho _{y}=\rho $. (This is
an easy consequence of Birkhoff ergodic theorem.)

\item There exists a set $\mathcal{J}\subset Y$ of total probability (that
is $\rho \left( \mathcal{J}\right) =1$ for every invariant measure $\rho $)
such that for all $y\in \mathcal{J}$ there exists an ergodic occupation
measure $\rho _{y}$.
\end{enumerate}

\begin{proposicao}
\label{proplyapreglyapint}Given $w\in \mathcal{W}$ there exists an invariant
ergodic measure $\mu ^{w}$ on $\mathbb{E}$ with $\pi _{*}\mu ^{w}=\nu $ such
that
\[
\int qd\mu ^{w}=w^{-1}H_{\mathrm{Ly}}
\]
and $\mu ^{w}(\mathrm{st}(w))=1$.
\end{proposicao}

\begin{profe}
If $\xi \in \mathrm{st}(w)$ then
\[
\lambda (\xi )=\lim_{k\rightarrow +\infty }\frac{1}{k}\mathsf{a}(k,\xi
)=w^{-1}H_{\mathrm{Ly}}
\]%
and since $\mathsf{a}(k,\xi )$ is a cocycle it follows that there exists an
occupation measure $\mu _{\xi }$ such that
\[
w^{-1}H_{\mathrm{Ly}}=\int qd\mu _{\xi }.
\]%
Note that $\pi _{\ast }(\mu _{\xi })$ is an occupation measure $\rho _{x}$
with $x=\pi (\xi )$. Since $\nu $ is ergodic, for $\nu $-almost all $x$, $%
\rho _{x}=\nu $ and hence we can choose $\xi $ with $\pi _{\ast }\left( \mu
_{\xi }\right) =\nu $.

It is not clear in advance that $\mu _{\xi }$ is ergodic. Nevertheless we
can decompose $\mu _{\xi }$ into ergodic components $\theta _{\eta }$ with $%
\eta $ ranging through a set $\mathcal{A}$ of $\mu _{\xi }$ total
probability, that is,
\[
\mu _{\xi }(\cdot )=\int \theta _{\eta }(\cdot )d\mu _{\xi }(\eta ).
\]
Since $\pi _{*}\mu _{\xi }=\nu $ it follows that $\pi _{*}\theta _{\eta
}=\nu $ for $\mu _{\xi }$-almost all $\eta $.

We claim that
\[
w^{-1}H_{\mathrm{Ly}}=\int qd\theta _{\eta }
\]%
for almost all $\eta \in \mathcal{A}$. In fact,
\[
w^{-1}H_{\mathrm{Ly}}=\int_{\mathbb{E}}\left( \int qd\theta _{\eta }\right)
d\mu _{\xi }(\eta ).
\]%
Hence $w^{-1}H_{\mathrm{Ly}}$ belongs to the convex closure of the set $%
\left\{ \int qd\theta _{\eta }\in \mathfrak{a};\,\eta \in \mathcal{A}%
\right\} $. However, by Proposition \ref{propmedinvconjest}, for any ergodic
$\theta _{\eta }$ there exists $u\in \mathcal{W}$ such that $\int qd\theta
_{\eta }=u^{-1}H_{\mathrm{Ly}}$, so that $w^{-1}H_{\mathrm{Ly}}$ is a convex
combination of points of the orbit $\mathcal{W}\cdot H_{\mathrm{Ly}}$. But
this is possible only if $\int qd\theta _{\eta }=w^{-1}H_{\mathrm{Ly}}$ for
almost all $\eta $, because the convex closure of the orbit $\mathcal{W}%
\cdot H_{\mathrm{Ly}}$ is a polyhedron whose vertices (extremal points) are
the points of the orbit. Hence there exists $\mu ^{w}$ yielding the Lyapunov
exponent $w^{-1}H_{\mathrm{Ly}}$. Finally, the equality $\mu ^{w}(\mathrm{st}%
(w))=1$ follows by the previous proposition.
\end{profe}

Now, we select two special kinds of ergodic measures on the flag bundles.

\begin{defi}
An ergodic measure $\mu $ on the maximal flag bundle $\mathbb{E}$ is said to
be an attractor measure for the flow if $\int qd\mu \in \mathrm{cl}\,%
\mathfrak{a}^{+}$. A measure $\mu _{\Theta }$ in $\mathbb{E}_{\Theta }$ is
an attractor measure if $\mu _{\Theta }={\pi _{\Theta }}_{\ast }\mu $ with $%
\mu $ attractor in $\mathbb{E}$, where $\pi _{\Theta }:\mathbb{E}\rightarrow
\mathbb{E}_{\Theta }$ is the canonical projection. \newline
Similarly, a measure $\mu $ in $\mathbb{E}$ is a repeller measure if $\int
qd\mu ^{w}\in -\mathrm{cl}\,\mathfrak{a}^{+}$, and $\mu _{\Theta }$ in $%
\mathbb{E}_{\Theta }$ is repeller if $\mu _{\Theta }={\pi _{\Theta }}_{\ast
}\mu $ with $\mu $ repeller in $\mathbb{E}$.
\end{defi}

Proposition \ref{proplyapreglyapint} ensures the existence of both attractor
and repeller measures.

\begin{proposicao}
\label{propatracrep}A repeller measure is an attractor measure for the
backward flow.
\end{proposicao}

\begin{profe}
Let $\mu $ be a repeller measure on $\mathbb{E}$ and write $q^{-}\left(
\cdot \right) =\mathsf{a}\left( -1,\cdot \right) $. Then by the cocycle
property $q^{-}(\xi )=-\mathsf{a}(1,\phi _{-1}(\xi ))=-q(\phi _{-1}(\xi ))$,
so that
\[
\int q^{-}d\mu =-\int qd\mu \in \mathrm{cl}\,\mathfrak{a}^{+}
\]
because $\int qd\mu \in -\mathrm{cl}\,\mathfrak{a}^{+}$. Thus $\mu $ is an
attractor measure for the backward flow. This proves the statement on the
maximal flag bundle $\mathbb{E}$. On the other bundles the result follows by
definition.
\end{profe}

Now we relate the supports of the attractor and repeller measures with the
decomposition given by the Multiplicative Ergodic Theorem on the flag
bundles $\mathbb{E}_{\Theta _{\mathrm{Ly}}\left( \nu \right) }$ and in its
dual $\mathbb{E}_{\Theta _{\mathrm{Ly}}^{\ast }\left( \nu \right) }$. This
decomposition is given by sections $\xi $ and $\xi ^{\ast }$ of $\mathbb{E}%
_{\Theta _{\mathrm{Ly}}\left( \nu \right) }$ and $\mathbb{E}_{\Theta _{%
\mathrm{Ly}}^{\ast }\left( \nu \right) }$, respectively.

We write simply $\Theta _{\mathrm{Ly}}=\Theta _{\mathrm{Ly}}\left( \nu
\right) $ and distinguish the several projections as: $\pi :\mathbb{E}%
\rightarrow X$, $\pi _{\Theta _{\mathrm{Ly}}}:\mathbb{E}\rightarrow \mathbb{E%
}_{\Theta _{\mathrm{Ly}}}$, $\pi _{\Theta _{\mathrm{Ly}}^{\ast }}:\mathbb{E}%
\rightarrow \mathbb{E}_{\Theta _{\mathrm{Ly}}^{\ast }}$ and $p$ for either $%
\mathbb{E}_{\Theta _{\mathrm{Ly}}}\rightarrow X$ or $\mathbb{E}_{\Theta _{%
\mathrm{Ly}}^{\ast }}\rightarrow X$.

Let $\mu $ be a repeller measure on $\mathbb{E}$ and put $\mu _{\Theta _{%
\mathrm{Ly}}^{*}}={\pi _{\Theta _{\mathrm{Ly}}^{*}}}_{*}\left( \mu \right) $
for the corresponding repeller measure on $\mathbb{E}_{\Theta _{\mathrm{Ly}%
}^{*}}$. We have $p_{*}(\mu _{\Theta _{\mathrm{Ly}}^{*}})=\nu $ because $%
p\circ \pi _{\Theta _{\mathrm{Ly}}^{*}}=\pi $ and $\pi _{*}\mu =\nu $. Hence
we can desintegrate $\mu _{\Theta _{\mathrm{Ly}}^{*}}$ with respect to $\nu $
to get
\[
\mu _{\Theta _{\mathrm{Ly}}^{*}}(\cdot )=\int_{X}\rho _{x}(\cdot )d\nu (x),
\]
where $x\in X\mapsto \rho _{x}\in \mathbb{M}^{+}\left( \mathbb{E}_{\Theta _{%
\mathrm{Ly}}^{*}}\right) $ is a measurable map into the space of probability
measures on $\mathbb{E}_{\Theta _{\mathrm{Ly}}^{*}}$.

\begin{lema}
For $\nu $-almost all $x\in X$ the component $\rho _{x}$ in the above
desintegration is a Dirac measure at $\xi ^{*}(x)$, that is, $\rho
_{x}=\delta _{\xi ^{*}\left( x\right) }$.
\end{lema}

\begin{profe}
Let $Z$ be the Borel set
\[
Z=\left\{ \mathrm{im}\,\xi ^{\ast }\right\} ^{c}=\mathbb{E}_{\Theta _{%
\mathrm{Ly}}^{\ast }}\backslash \left\{ \mathrm{im}\,\xi ^{\ast }\right\} .
\]%
Then
\[
\mu _{\Theta _{\mathrm{Ly}}^{\ast }}(Z)=\mu (\pi _{\Theta _{\mathrm{Ly}%
}^{\ast }}^{-1}(Z))=\mu (\mathbb{E}\backslash \mathrm{st}(w_{0}))=0,
\]%
because $\mu $ is a repeller measure. However,
\[
0=\mu _{\Theta _{\mathrm{Ly}}^{\ast }}(Z)=\int_{X}\rho _{x}(Z)d\nu (x),
\]%
and since $\rho _{x}$ is supported at $\pi ^{-1}(x)$, it follows that $\rho
_{x}(\mathbb{E}_{\Theta _{\mathrm{Ly}}^{\ast }}\setminus \{\xi ^{\ast
}(x)\})=0$, for $\nu $-almost all $x\in X$.
\end{profe}

This lemma shows also that a repeller measure on the dual flag manifold $%
\mathbb{E}_{\Theta _{\mathrm{Ly}}^{*}}$ is unique. Now we can apply the same
argument for the reverse flow $\phi _{-t}$, and get a similar result now for
an attracting measure on $\mathbb{E}_{\Theta _{\mathrm{Ly}}}$ with $\xi ^{*}$
replaced by $\xi $.

For later reference we summarize these facts in the following proposition.

\begin{proposicao}
\label{propmedatr}There exists a unique attractor measure $\mu _{\Theta _{%
\mathrm{Ly}}}^{+}$ for $\phi _{t}$ in its flag type $\mathbb{E}_{\Theta _{%
\mathrm{Ly}}}$, which is a Dirac measure on $\xi \left( x\right) $, that is,
it desintegrates as
\[
\mu _{\Theta _{\mathrm{Ly}}}^{+}(\cdot )=\int \delta _{\xi (x)}(\cdot )d\nu
(x)
\]
with respect to $\nu $. There exists also a unique repeller measure $\mu
_{\Theta _{\mathrm{Ly}}}^{-}$ on $\mathbb{E}_{\Theta _{\mathrm{Ly}}^{*}}$
which is Dirac at $\xi ^{*}$.
\end{proposicao}

\begin{corolario}
\label{coruniquemedatr}There exists a unique attractor (respectively
repeller) measure in $\mathbb{E}_{\Theta }$ if $\Theta _{\mathrm{Ly}}\subset
\Theta $ (respectively $\Theta _{\mathrm{Ly}}^{*}\subset \Theta $).
\end{corolario}

\begin{profe}
This is because the projection $\mathbb{E}\rightarrow \mathbb{E}_{\Theta }$
factors through $\mathbb{E}_{\Theta _{\mathrm{Ly}}}$ if $\Theta _{\mathrm{Ly}%
}\subset \Theta $: $\mathbb{E}\rightarrow \mathbb{E}_{\Theta _{\mathrm{Ly}%
}}\rightarrow \mathbb{E}_{\Theta }$. Hence a measure in $\mathbb{E}_{\Theta
} $ is attractor if and only if it is the projection of the attractor
measure in $\mathbb{E}_{\Theta _{\mathrm{Ly}}}$.
\end{profe}

\section{Morse decomposition $\times $ Oseledets decomposition\label%
{secversusdec}}

In this section we use the concepts of attractor and repeller measures
developed above to relate the Oseledets decomposition and the Morse
decomposition on a flag bundle $\mathbb{E}_{\Theta }$, as well as the
Lyapunov spectrum \ $H_{\mathrm{Ly}}$ and the Morse spectrum $\Lambda _{%
\mathrm{Mo}}$.

First we have the following consequence of Proposition \ref%
{proplyapreglyapint}.

\begin{proposicao}
\label{proplmbdamocontido}Suppose that $\alpha \left( \Lambda _{\mathrm{Mo}%
}\right) =0$ for all $\alpha \in \Theta _{\mathrm{Mo}}$. Then $\Theta _{%
\mathrm{Ly}}\left( \rho \right) =\Theta _{\mathrm{Mo}}$ for every ergodic
measure $\rho $ on the base space.
\end{proposicao}

\begin{profe}
As checked in Section \ref{seclyapmorse} we have $\Theta _{\mathrm{Ly}%
}\left( \rho \right) \subset \Theta _{\mathrm{Mo}}$. On the other hand by
Proposition \ref{proplyapreglyapint} any regular Lyapunov exponent is a
Morse exponent, that is, $H_{\mathrm{Ly}}\left( \rho \right) \subset \Lambda
_{\mathrm{Mo}}$. So that $\alpha \left( H_{\mathrm{Ly}}\left( \rho \right)
\right) =0$ if $\alpha \in \Theta _{\mathrm{Mo}}$, showing that $\Theta _{%
\mathrm{Mo}}\subset \Theta _{\mathrm{Ly}}\left( \rho \right) $.
\end{profe}

Now we look at the decompositions of the flag bundles.

\begin{proposicao}
Let $\mu $ be an attractor measure on $\mathbb{E}$. Then its support $%
\mathrm{supp}\mu $ is contained in the unique attractor component $\mathcal{M%
}^{+}$ of the finest Morse decomposition.
\end{proposicao}

\begin{profe}
Each point in $\mathrm{supp}\mu $ is recurrent and hence belongs to the set
of chain recurrent points which is the union of the Morse components.

Now since $\mu $ is an attractor measure, by definition the integral
\[
\lambda _{\mu }=\int qd\mu \in \mathrm{cl}\mathfrak{a}^{+}.
\]%
This integral is the $\mathfrak{a}$-Lyapunov exponent of $\mu $-almost all $%
z\in \mathrm{supp}\mu $. Hence is contained in the Morse spectrum. Actually,
$\lambda _{\mu }\in \Lambda _{\mathrm{Mo}}\left( \mathcal{M}^{+}\right) $,
the Morse spectrum of $\mathcal{M}^{+}$, because this is the only Morse
component whose spectrum meets $\mathrm{cl}\mathfrak{a}^{+}$. Therefore for $%
\mu $-almost all $z\in \mathrm{supp}\mu $, $z\in \mathcal{M}^{+}$. Since $%
\mathcal{M}^{+}$ is compact it follows that $\mathrm{supp}\mu \subset
\mathcal{M}^{+}$.
\end{profe}

By taking the backward flow we get a similar result for the repeller
measures.

\begin{proposicao}
Let $\mu $ be a repeller measure on $\mathbb{E}$. Then its support $\mathrm{%
supp}\mu $ is contained in the unique repeller component $\mathcal{M}^{-}$
of the finest Morse decomposition.
\end{proposicao}

\begin{proposicao}
\label{prooselcontmor}Let $\mathcal{O}$ be a component of the Oseledets
decomposition in a flag bundle $\mathbb{E}_{\Theta }$. Then there exists a
component $\mathcal{M}$ of the Morse decomposition of $\mathbb{E}_{\Theta }$
such that $\mathcal{O}\subset \mathcal{M}$.
\end{proposicao}

\begin{profe}
Let $\mu _{\Theta _{\mathrm{Ly}}}$ the only attractor measure in $\mathbb{E}%
_{\Theta _{\mathrm{Ly}}}$ and $\mu _{\Theta _{\mathrm{Ly}}^{\ast }}$ the
repeller measure in $\mathbb{E}_{\Theta _{\mathrm{Ly}}^{\ast }}$. These are
projections of attractor and repeller measures on $\mathbb{E}$. Hence the
above lemmas imply that $\mathrm{supp}\mu _{\Theta _{\mathrm{Ly}}}\subset
\mathcal{M}_{\Theta _{\mathrm{Ly}}}^{+}$ and $\mathrm{supp}\mu _{\Theta _{%
\mathrm{Ly}}^{\ast }}\subset \mathcal{M}_{\Theta _{\mathrm{Ly}}^{\ast }}^{+}$%
. However we checked in Section \ref{seclyapmorse} that $\Theta _{\mathrm{Ly}%
}\subset \Theta _{\mathrm{Mo}}$. Hence by Lemma \ref{lemfixinclu}, we
conclude that the fixed point set -- in any flag bundle -- of the section
Oseledet section $\chi _{\mathrm{Ly}}$ is contained in the fixed point set
of the Morse section $\chi _{\mathrm{Mo}}$. This means that the Oseletet
components are contained in the Morse components.
\end{profe}

\vspace{12pt}%

\noindent \textbf{Remark:} It is proved in \cite{ck}, Corollary 5.5.17, that
the Oseledets decomposition is contained in the Morse decomposition for a
linear flow on a vector bundle.

\section{Lyapunov exponents in the tangent bundle $T^{f}\mathbb{E}_{\Theta _{%
\mathrm{Ly}}}$}

\label{lyextgbun}

A fiber of a flag bundle $\mathbb{E}_{\Theta }$ is differentiable manifold
and hence has a tangent bundle. Gluing together the tangent bundles to the
fibers of $\mathbb{E}_{\Theta }$ we get a vector bundle $T^{f}\mathbb{E}%
_{\Theta }\rightarrow \mathbb{E}_{\Theta }$ over $\mathbb{E}_{\Theta }$. The
flow $\phi _{t}$ on $\mathbb{E}_{\Theta }$ is differentiable along the
fibers with differential map $\psi _{t}$, a linear map of the vector bundle $%
T^{f}\mathbb{E}_{\Theta }$. (See \cite{conleyflags} for a construction of
this vector bundle as an associated bundle $Q\times _{G}V$.)

We look here at the Lyapunov exponents for the linear flow $\psi _{t}$ on $%
T^{f}\mathbb{E}_{\Theta _{\mathrm{Ly}}}$ with respect to an attractor
measure of $\phi _{t}$. These Lyapunov exponents will be used in the proof
of the main technical lemma to describe the $\omega $-limit sets (w.r.t. $%
\phi _{t}$) in the flag bundles.

Equip the bundle $T^{f}\mathbb{E}_{\Theta }$ with a Riemannian metric $%
\langle \cdot ,\cdot \rangle $, which can be built from a $K$-reduction $R$
of the principal bundle $Q$. (Roughly, the metric $\langle \cdot ,\cdot
\rangle $ is constructed by piecing together $K$-invariant metrics on the
fibers. See \cite{conleyflags} for details.)

\begin{proposicao}
Let $\nu $ be an ergodic measure on $X$ and denote by $\mu $ its attractor
measure on the bundle $\mathbb{E}_{\Theta _{\mathrm{Ly}}}$. Let $H\left( \nu
\right) \in \mathrm{cl}\mathfrak{a}^{+}$ be the polar exponent of $\nu $.
Then the Lyapunov spectrum of $\psi $ w.r.t. $\mu $ is $\mathrm{ad}\left(
H\left( \mu \right) \right) _{\left\vert \mathfrak{n}_{\Theta _{\mathrm{Ly}%
}}^{-}\right. }$, which is a diagonal linear map of $\mathfrak{n}_{\Theta _{%
\mathrm{Ly}}}^{-}$ (that is an element of a Weyl chamber $\mathrm{cl}{%
\mathfrak{a}}^{+}$ of $\mathfrak{gl}\left( \mathfrak{n}_{\Theta _{\mathrm{Ly}%
}}^{-}\right) $).
\end{proposicao}

\begin{profe}
Denote by $\mathcal{O}\left( \nu \right) $ the $Z_{H_{\mathrm{Ly}}}$%
-measurable reduction of $Q$, corresponding to the Oseledets section of $\nu
$. This reduction is a principal bundle with structural group $Z_{H_{\mathrm{%
Ly}}}$ over a set $\Omega \subset X$ with $\nu \left( \Omega \right) =1$
(see \cite{alvsm}).

The section $\xi _{\mathrm{Ly}}:\Omega \rightarrow \mathbb{E}_{\Theta _{%
\mathrm{Ly}}}$ gives a desintegration of $\mu $ with respect to $\nu $ by
Dirac measures. Let $\Omega ^{\#}$ be the image of this section. Then the
restriction of $T^{f}\mathbb{E}_{\Theta _{\mathrm{Ly}}}$ to $\Omega ^{\#}$
is a vector bundle $T^{f}\Omega ^{\#}\rightarrow \Omega ^{\#}$, which is
invariant by the differential flow $\psi _{t}$.

We can build the vector bundle $T^{f}\Omega ^{\#}$ as an associated bundle $%
\mathcal{O}\left( \nu \right) $ through the adjoint representation $\theta $
of $Z_{H_{\mathrm{Ly}}}$ in $\mathfrak{n}_{\Theta _{\mathrm{Ly}}}^{-}$.

Then if we take compatible Cartan decompositions of $\theta \left( Z_{H_{%
\mathrm{Ly}}}\right) $ and $\mathrm{Gl}\left( \mathfrak{n}_{\Theta _{\mathrm{%
Ly}}}^{-}\right) $ it follows that the polar exponent of $\psi _{t}$ is
precisely $\theta \left( H\left( \mu \right) \right) $. By the constructions
of Section 8 of \cite{alvsm} it follows that the Lyapunov exponents of $\psi
_{t}$ are the eigenvalues of $\theta \left( H\left( \mu \right) \right) $ as
linear maps of $\mathfrak{n}_{\Theta _{\mathrm{Ly}}}^{-}$.
\end{profe}

\begin{corolario}
\label{corlyapvertneg}Suppose $\Theta _{\mathrm{Ly}}\subset \Theta $. Then
the Lyapunov exponents of $\psi _{t}$ in $T^{f}\mathbb{E}_{\Theta }$ with
respect to the attractor measure $\mu $ are strictly negative.
\end{corolario}

Later on we will combine this corollary with the following general fact
about Lyapunov exponents on vector bundles. Let $p:V\rightarrow X$ be a
continuous vector bundle endowed with a norm $\left\vert \left\vert \cdot
\right\vert \right\vert $. Let $\Phi _{n}$ be a continuous linear flow on $V$%
. If $\nu $ is a $\Phi $-invariant ergodic measure on the base $X$ then $%
\Phi $ has a Lyapunov spectrum $H_{\mathrm{Ly}}\left( \nu \right) =\{\lambda
_{1}\geq \cdots \geq \lambda _{n}\}$ with respect to $\nu $, as ensured by
the Multiplicative Ergodic Theorem. The following lemma may be well known.
For the sake of completeness we prove it here using the Morse spectrum of
the linear flow.

\begin{lema}
\label{lemfibradovetorial}Supppose that for every $\Phi $-invariant ergodic
measure $\nu $ on $X$ the spectrum with respect to $\nu $ is strictly
negative. Then for every $v\in V$,
\[
\lim_{n\rightarrow +\infty }\left| \left| \Phi _{n}v\right| \right| =0.
\]
\end{lema}

\begin{profe}
Let $\mathbb{P}V\rightarrow X$ be the projective bundle of $V$. The cocycle $%
\rho \left( n,v\right) =\frac{\left\vert \left\vert \Phi _{n}v\right\vert
\right\vert }{\left\vert \left\vert v\right\vert \right\vert }$ on $V$
induces the additive cocycle $a\left( n,\eta \right) =\log \rho \left(
n,\eta \right) $, $\eta \in \mathbb{P}V$, whose asymptotics give the
Lyapunov spectrum of $\Phi $. Write $q\left( \cdot \right) =a\left( 1,\cdot
\right) $. Then by general results on Morse spectrum of an additive cocycle
(see \cite{smsec}, Section 3, and references therein), the Morse spectrum of
$a$ is a union of intervals whose extreme points are integrals $\int q\left(
x\right) \mu \left( dx\right) $ with respect to ergodic invariant measures $%
\mu $ for the flow on $\mathbb{P}V$. By the Birkhoff ergodic theorem it
follows that for $\mu $-almost all $\eta \in \mathbb{P}V$,
\[
\lim_{n\rightarrow \infty }\frac{1}{n}a\left( n,\eta \right)
=\lim_{n\rightarrow \infty }\frac{1}{n}\sum_{k=0}^{n-1}q\left( \Phi _{k}\eta
\right) =\int q\left( z\right) \mu \left( dz\right) .
\]%
On the other hand the projection $p_{\ast }\mu =\nu $ is ergodic on the base
$X$. Hence by assumption the spectrum with respect to $\nu $, given by the
multiplicative ergodic theorem is strictly negative. This means that for $%
\nu $-allmost all $x\in X$, $\lim_{n\rightarrow \infty }\frac{1}{n}a\left(
n,\eta \right) $ exists for every $\eta \in p^{-1}\{x\}$ and is strictly
negative. Combining these two facts we conclude that
\[
\int q\left( z\right) \mu \left( dz\right) <0,
\]%
and therefore  the Morse spectrum is contained in $\left( -\infty ,0\right) $%
.

Now, for every $\eta \in \mathbb{P}V$, $\lim \sup_{n\rightarrow +\infty }%
\frac{1}{n}a\left( n,\eta \right) $ belongs to the Morse spectrum (see \cite%
{ck}, Theorem 5.3.6). Hence for every $0\neq v\in V$,
\[
\lim \sup_{n\rightarrow +\infty }\frac{1}{n}\log \left\vert \left\vert \Phi
_{n}v\right\vert \right\vert <0.
\]%
This implies that for large $n$, $\left\vert \left\vert \Phi
_{n}v\right\vert \right\vert <e^{cn}$, $c<0$, proving the lemma.
\end{profe}

Applying the lemma for the backward flow we have

\begin{corolario}
\label{corfibradovetorialneg}With the same assumptions of the lemma we have $%
\lim_{n\rightarrow -\infty }\left\vert \left\vert \Phi _{n}v\right\vert
\right\vert =\infty $ se $v\neq 0$.
\end{corolario}

\section{Main technical lemma}

\label{teclemma}

\begin{lema}
\label{lemmaintech}Let $\mathbb{E}_{\Theta }$ be a flag manifold with dual $%
\mathbb{E}_{\Theta ^{\ast }}$. Suppose there are  three compact $\phi $%
-invariant subsets $A,B\subset \mathbb{E}_{\Theta }$ and $C\subset \mathbb{E}%
_{\Theta ^{\ast }}$ that project onto $X$ and such that

\begin{enumerate}
\item $A\cap B=\emptyset $.

\item $B^{c}$ is the set of elements transversal to $C$. (That is an element
$v\in \mathbb{E}_{\Theta }$ belongs to $B$ if and only if it is not
transversal to some $w\in C$ in the same fiber as $v$.)

\item For any ergodic measure $\rho $ for the flow on the base space $X$ we
have $\Theta _{\mathrm{Ly}}\left( \rho \right) \subset \Theta $. By
Corollary \ref{coruniquemedatr}, this implies that there is a unique
attractor measure $\mu _{\Theta }^{+}\left( \rho \right) $ for $\rho $ on $%
\mathbb{E}_{\Theta }$ and a unique repeller measure $\mu _{\Theta ^{\ast
}}^{-}\left( \rho \right) $ on $\mathbb{E}_{\Theta ^{\ast }}$.

\item For any ergodic measure $\rho $ on $X$, $\mathrm{supp}\mu _{\Theta
}^{+}\left( \rho \right) \subset A$ and $\mathrm{supp}\mu _{\Theta ^{\ast
}}^{-}\left( \rho \right) \subset C$.
\end{enumerate}

Then $\left( A,B\right) $ is an attractor-repeller pair on $\mathbb{E}%
_{\Theta }$. That is, the $\omega $-limit $\omega \left( v\right) \subset A$
if $v\notin B$ and $\omega ^{\ast }\left( v\right) \subset B$ if $v\notin A$.
\end{lema}

The proof of this lemma will be done in several steps. Before starting we
define a fourth set $D\subset Q\times _{G}\left( \mathbb{F}_{\Theta }\times
\mathbb{F}_{\Theta ^{\ast }}\right) $ by
\[
D=\pi _{1}^{-1}\left( A\right) \cap \pi _{2}^{-1}\left( C\right)
\]%
where $\pi _{1}:Q\times _{G}\left( \mathbb{F}_{\Theta }\times \mathbb{F}%
_{\Theta ^{\ast }}\right) \rightarrow \mathbb{E}_{\Theta }$ and $\pi
_{2}:Q\times _{G}\left( \mathbb{F}_{\Theta }\times \mathbb{F}_{\Theta ^{\ast
}}\right) \rightarrow \mathbb{E}_{\Theta ^{\ast }}$ are the projections.
This set is compact and invariant, and by the transversality given by the
first and second conditions in the lemma we can view $D$ as a compact subset
of the bundle
\[
\mathcal{A}_{\Theta }=Q\times _{G}\mathrm{Ad}\left( G\right) H_{\Theta }
\]%
where $\Theta =\{\alpha \in \Sigma :\alpha \left( H_{\Theta }\right) =0\}$.

Now, to start the proof fix $x\in X$ that has a periodic orbit $\mathcal{O}%
\left( x\right) $. Then the Oseledets decomposition coincides with the Morse
decomposition above $\mathcal{O}\left( x\right) $, which by Section \ref%
{secperflow} is built from the dynamics of the action of a $g_{x}\in G$.
Clearly the homogeneous measure $\theta $ on the periodic orbit is an
ergodic invariant measure on $X$. By the third condition of the lemma $g_{x}
$ has one attractor fixed point at $\mathbb{F}_{\Theta }$, say $b^{+}$, and
a repeller fixed point $b^{-}\in \mathbb{F}_{\Theta ^{\ast }}$. The Morse
decomposition of $g_{x}$ is the union of $\{b^{+}\}$ with subsets whose
elements are not transversal to $b^{-}$. It follows that the attractor $\mu
_{\Theta }^{+}\left( \theta \right) $ and the repeller measures $\mu
_{\Theta ^{\ast }}^{-}\left( \theta \right) $ are homogeneous measures on
periodic orbits. By the fourth condition of the lemma $\mathrm{supp}\mu
_{\Theta }^{+}\left( \theta \right) \subset A$ and $\mathrm{supp}\mu
_{\Theta ^{\ast }}^{-}\left( \theta \right) \subset C$, so that by the
second condition the Morse decomposition above $\mathcal{O}\left( x\right) $
is the union of $\mathrm{supp}\mu _{\Theta }^{+}\left( \theta \right) $ with
subsets contained in $B$. Hence, if $v$ is in the fiber above $x$ then $%
\omega \left( v\right) \subset A$ if $v\notin B$ while $\omega ^{\ast
}\left( v\right) \subset B$ if $v\notin A$.

Therefore, from now on we look at $\omega $-limits $\omega \left( v\right) $
and $\omega ^{\ast }\left( v\right) $ assuming that the orbit $\mathcal{O}%
\left( x\right) $ of $x=\pi \left( v\right) $ is not periodic, that is, the
map $n\in \mathbb{Z}\mapsto x_{n}=\phi _{n}\left( x\right) \in X$ is
injective and $\mathcal{O}\left( x\right) $ is in bijection with $\mathbb{Z}$%
.

To prove that the $\omega $-limits are contained in $A$ we will use the
following consequence of Lemma \ref{lemfibradovetorial}.

\begin{lema}
Given $v\in A$ let $w\in T_{v}^{f}\mathbb{E}_{\Theta _{\mathrm{Ly}}}$ be a
tangent vector at $v$. Then $\lim_{t\rightarrow +\infty }\left| \left| \psi
_{t}w\right| \right| =0$.
\end{lema}

\begin{profe}
Is an immediate consequence of Lemma \ref{lemfibradovetorial}, combined with
the third and fourth conditions of the Lemma.
\end{profe}

Now, above the nonperiodic orbit $\mathcal{O}\left( x\right) $ we reduce the
flow to just a sequence $g_{n}$ of elements of the subgroup $Z_{H_{\Theta }}$%
. The construction is the following: Start with an element $\eta _{0}\in D$
in the fiber over $x$. The orbit $\mathcal{O}\left( \eta \right) $ is the
sequence $\eta _{n}=\phi _{n}\left( \eta _{0}\right) $, $n\in \mathbb{Z}$,
that can be viewed as a section over $\mathcal{O}\left( x\right) $ by $%
x_{n}\mapsto \eta _{n}$. The elements of the associated bundle $\mathcal{A}%
_{\Theta }$ are written as $p\cdot H_{\Theta }$, $p\in R$, where as before $%
R $ is the $K$-reduction of $Q\rightarrow X$. Hence there exists a sequence $%
p_{n}\in R$ such that $\eta _{n}=p_{n}\cdot H_{\Theta }$.

Since $p_{n+m}$ and $\phi _{n}\left( p_{m}\right) $ are in the same fiber,
we have $\phi _{n}\left( p_{m}\right) =p_{n+m}\cdot g_{n,m}$ with $%
g_{n,m}\in G$, $n,m\in \mathbb{Z}$. Actually, $g_{n,m}\in Z_{H_{\Theta }}$
because $\phi _{n}\left( \eta _{m}\right) =\eta _{n+m}$, so that
\[
p_{n+m}\cdot \mathrm{Ad}\left( g_{n,m}\right) H_{\Theta }=\phi _{n}\left(
\eta _{m}\right) =\eta _{n+m}=p_{n+m}\cdot H_{\Theta }.
\]

We write $\xi _{n}$ and $\xi _{n}^{*}$ for the projections of $\eta _{n}$
into $\mathbb{E}_{\Theta }$ and $\mathbb{E}_{\Theta ^{*}}$, respectively. By
definition of $D$ we have $\xi _{n}\in A$ and $\xi _{n}^{*}\in C$. Hence by
the second assumption of the lemma the elements in $\mathbb{E}_{\Theta }$
that are not transversal to $\xi _{n}^{*}$ are contained in $B$. In other
words

\begin{lema}
Take $v\notin B$ in the fiber of $x$. Then $v$ is transversal to $\xi
_{0}^{*}$.
\end{lema}

Now we use Lyapunov exponents of the lifting $\psi _{n}$ of $\phi _{n}$ to $%
T^{f}\mathbb{E}_{\Theta }$ to show that $\omega \left( v\right) \subset A$
if $v\notin B$ is in the fiber of $x$.

To do that we first note that if the starting element $\eta _{0}\in \mathcal{%
A}_{\Theta }$ is written as $\eta _{0}=p\cdot H_{\Theta }$, $p\in R$, then
the set of points that are transversal to $\xi _{0}^{*}$ is given
algebraically by
\[
T=p\cdot \left( N_{\Theta }^{-}\cdot b_{0}\right) =\{p\cdot nb_{0}:n\in
N_{\Theta }^{-}\}
\]
where $b_{0}$ is the origin of the flag $\mathbb{F}_{\Theta }$ and $%
N_{\Theta }^{-}$ is the nilpotent subgroup with Lie algebra $\mathfrak{n}%
_{\Theta }^{-}=\sum_{\alpha \notin \langle \Theta \rangle ,\alpha <0}%
\mathfrak{g}_{\alpha }$ (lower triangular matrices).

Since $\exp :\mathfrak{n}_{\Theta }^{-}\rightarrow N_{\Theta }^{-}$ is a
diffeomorphism we have also $T=\{p\cdot \left( \exp Y\cdot b_{0}\right)
:Y\in \mathfrak{n}_{\Theta }^{-}\}$.

The action of $\phi _{n}$ on an element $p\cdot \left( \exp Y\cdot
b_{0}\right) \in T$ is given as follows: Put $g_{n}=g_{n,0}\in Z_{H_{\Theta
}}$. Then, as remarked above, $\phi _{n}\left( p\right) =p_{n}\cdot g_{n}$,
so that
\[
\phi _{n}\left( p\cdot \left( \exp Y\cdot b_{0}\right) \right) =p_{n}\cdot
\left( g_{n}\exp Yg_{n}^{-1}\cdot g_{n}b_{0}\right) .
\]
But $g_{n}b_{0}=b_{0}$ because $g_{n}\in Z_{H}$, and $g_{n}\exp
Yg_{n}^{-1}=\exp \left( \mathrm{Ad}\left( g_{n}\right) Y\right) $. Hence
\begin{equation}
\phi _{n}\left( p\cdot \left( \exp Y\cdot b_{0}\right) \right) =p_{n}\cdot
\exp \left( \mathrm{Ad}\left( g_{n}\right) Y\right) b_{0}.  \label{forfinexp}
\end{equation}

The next lemma relates this action with the lifting $\psi _{n}$ of $\phi
_{n} $ to the tangent space $T^{f}\mathbb{E}_{\Theta }$.

\begin{lema}
\label{lemlimzerotangfibr} Given $Y\in \mathfrak{g}$, denote by $p \cdot Y$
the vertical tangent vector \linebreak $\frac{d}{dt}\left( p\cdot \left(
\exp tY\cdot b_{0}\right) \right) _{t=0}\in T_{p\cdot b_{0}}^{f}\mathbb{E}%
_{\Theta }$. Then $p\cdot Y$, $Y\in \mathfrak{n}_{\Theta }^{-}$, fulfill the
vertical tangent space $T_{p\cdot b_{0}}^{f}\mathbb{E}_{\Theta }$, and the
derivative $\psi _{n}$ of $\phi _{n}$ at $p\cdot b_{0}$ satisfies
\[
\psi _{n}\left( p\cdot Y\right) =p_{n}\cdot \mathrm{Ad}\left( g_{n}\right)
Y.
\]
\end{lema}

\begin{profe}
The fact that any tangent vector in $T_{p\cdot b_{0}}^{f}\mathbb{E}_{\Theta }
$ is given by $p\cdot Y$ for some $Y\in \mathfrak{n}_{\Theta }^{-}$ is due
to the fact that $N_{\Theta }^{-}\cdot b_{0}$ is an open submanifold of $%
\mathbb{F}_{\Theta }$. For the last statement we have
\begin{eqnarray*}
\psi _{n}\left( p\cdot Y\right)  &=&\frac{d}{dt}\phi _{n}\left( p\cdot
\left( \exp tY\cdot b_{0}\right) \right) _{t=0} \\
&=&p_{n}\cdot \frac{d}{dt}\left( \exp \left( t\mathrm{Ad}\left( g_{n}\right)
Y\right) b_{0}\right) _{t=0}=p_{n}\cdot \mathrm{Ad}\left( g_{n}\right) Y.
\end{eqnarray*}
\end{profe}

We are now prepared to prove that $\omega \left( v\right) \subset A$ if $%
v\notin B$ is in the fiber of $x$. We have $v=p\cdot \left( \exp Y\cdot
b_{0}\right) $ for some $Y\in \mathfrak{n}_{\Theta }^{-}$, so that by (\ref%
{forfinexp}) $\phi _{n}\left( v\right) =p_{n}\cdot \exp \left( \mathrm{Ad}%
\left( g_{n}\right) Y\right) b_{0}$.

Now by Corollary \ref{corlyapvertneg} and \ Lemma \ref{lemfibradovetorial}
we have $\lim_{n\rightarrow +\infty }\left\vert \left\vert \psi
_{n}w\right\vert \right\vert =0$ if $w\in T_{p\cdot b_{0}}^{f}\mathbb{E}%
_{\Theta }$. Taking $w=p\cdot Y$ we have $\left\vert \left\vert p\cdot
Y\right\vert \right\vert =\left\vert \left\vert Y\right\vert \right\vert $
because $p\in R$ and hence is an isometry between the flag manifold $\mathbb{%
F}_{\Theta }$ and the corresponding fiber of $\mathbb{E}_{\Theta }$ (see
\cite{conleyflags} for the construction of the norm in $T^{f}\mathbb{E}%
_{\Theta }$). Since the same remark holds for $p_{n}\in R$ we have $%
\left\vert \left\vert \psi _{n}w\right\vert \right\vert =\left\vert
\left\vert \mathrm{Ad}\left( g_{n}\right) Y\right\vert \right\vert $, so
that
\[
\mathrm{Ad}\left( g_{n}\right) Y\rightarrow 0\qquad \mathrm{and}\qquad \exp
\mathrm{Ad}\left( g_{n}\right) Y\rightarrow 1.
\]

This implies that if $d$ is the metric on $\mathbb{E}_{\Theta _{\mathrm{Ly}%
}} $ then $d\left( \phi _{n}\left( v\right) ,p_{n}\cdot b_{0}\right)
\rightarrow 0$. But $p_{n}\cdot b_{0}=\xi _{n} \in A$ as well as its limit
points, by invariance and compactness of $A$. Therefore we conclude that $%
\omega \left( v\right) \subset A$, if $v\notin B$.

We turn now to the proof that $\omega ^{*}\left( v\right) \subset B$ if $%
v\notin A$. Again with $v$ above a nonperiodic orbit.

Take a sequence $n_{k}\rightarrow -\infty $ such that $\phi _{n_{k}}v$
converges in $\mathbb{E}_{\Theta }$. Taking subsequences we assume the
convergences $p_{n_{k}}\rightarrow p\in R$, $\eta _{n_{k}}\rightarrow \eta $%
, $\xi _{n_{k}}\rightarrow \xi $ and $\xi _{n_{k}}^{*}\rightarrow \xi ^{*}$.

By invariance it is enough to take $v\notin B$, so that we can write $%
v=p_{0}\cdot \left( \exp Y\right) b_{0}$ with $Y\in \mathfrak{n}_{\Theta
}^{-}$ and $Y\neq 0$ (because $v\notin A$).

Taking subsequences again we assume that $g_{n_{k}}$ $\left( \exp Y\right)
b_{0}$ converges to $b_{1}\in \mathbb{F}_{\Theta }$. Since $g_{n}\in
Z_{H_{\Theta }}$, we have $g_{n}b_{0}=b_{0}$ and hence
\[
\left( \exp \mathrm{Ad}\left( g_{n_{k}}\right) Y\right)
b_{0}=g_{n_{k}}\left( \exp Y\right) b_{0}\rightarrow b_{1}.
\]

Now $\mathrm{Ad}\left( g_{n}\right) Y\rightarrow \infty $ in $\mathfrak{n}%
_{\Theta _{\mathrm{Ly}}}^{-}$, because the Lyapunov exponents for the
backward flow are $>0$. This implies that $b_{1}$ is not transversal to the
origin $b_{0}^{*}$ of $\mathbb{F}_{\Theta ^{*}}$. Hence, $p_{n_{k}}\cdot
b_{1}$ is not transversal to $p_{n_{k}}\cdot b_{0}^{*}$, so that $%
p_{n_{k}}\cdot b_{1}\in B$. But
\[
\phi _{n_k}v=p_{n_{k}}g_{n_{k}}\cdot \left( \exp Y\right)
b_{0}=p_{n_{k}}\cdot \left( \exp \mathrm{Ad}\left( g_{n_{k}}\right) Y\right)
b_{0}
\]
so that $\lim \phi _{n_k}v=p\cdot b_{1}$, showing that $\omega ^{*}\left(
v\right) \subset B$.

In conclusion we have compact invariant sets $A$ and $B$ that satisfy $%
\omega \left( v\right) \subset A$ and $\omega ^{*}\left( v\right) \subset B$
if $v\notin A\cup B$. Hence $A$ and $B$ define a Morse decomposition of $%
\mathbb{E}_{\Theta _{\mathrm{Ly}}}$ with $A$ the attractor component.

\section{Three conditions}

\label{threecon}

In this section we state three conditions that together are necessary and
sufficient to have equality between the Lyapunov and Morse decompositions
over an ergodic invariant measure.

Thus as in Section \ref{seclyapmorse}, let $\phi _{n}$ be a continuous
flow on a continuous principal bundle $\pi :Q\rightarrow X$ whose the
structural group $G$ is reductive and noncompact. We fix once and for all an
ergodic invariant measure $\nu $ on the base space having support $\mathrm{%
supp}\nu =X$. Then the  $\mathfrak{a}$-Lyapunov exponents of $\phi _{n}$
select a flag type, which is expressed by a subset $\Theta _{\mathrm{Ly}}$
of simple roots. The flow on $X$ is chain transitive so it also has a flag
type $\Theta _{\mathrm{Mo}}$ coming from the Morse decompositions and $%
\mathfrak{a}$-Morse spectrum.

We start by writting down the three conditions and check that they are
necessary. In the next section we prove that together they are also
sufficient to have $\Theta _{\mathrm{Ly}}=\Theta _{\mathrm{Mo}}$.

\subsection{Bounded section\label{cond1}}

The Oseledets' section is a measurable section $\chi _{\mathrm{Ly}}:\Omega
\rightarrow Q\times _{G}\mathcal{O}_{\mathrm{Ly}}$ of the associated bundle $%
Q\times _{G}\mathcal{O}_{\mathrm{Ly}}\rightarrow X$ above the set of full $%
\nu $-measure $\Omega $. The fiber of this bundle is the adjoint orbit $%
\mathcal{O}_{\mathrm{Ly}}=\mathrm{Ad}\left( G\right) H_{\mathrm{Ly}}$. This
section can be seen as an equivariant map $f_{\mathrm{Ly}}:Q_{\Omega
}\rightarrow \mathcal{O}_{\mathrm{Ly}}$ defined above $\Omega $, where $%
Q_{\Omega }=\pi ^{-1}\left( \Omega \right) $ and $\pi :Q\rightarrow X$ is
the projection.

Let $R\subset Q\rightarrow X$ be a (continuous) $K$-reduction of $Q$. Then
we say that the Oseledets' section is \textit{bounded} if

\begin{itemize}
\item $f_{\mathrm{Ly}}$ is bounded in $R_{\Omega }$.
\end{itemize}

This definition does not depend on the specific $K$-reduction because the
base space $X$ is assumed to be compact.

If $\Theta _{\mathrm{Ly}}=\Theta _{\mathrm{Mo}}$ then we can take $H_{%
\mathrm{Mo}}=H _{\mathrm{Ly}}$ and $\chi _{\mathrm{Mo}}=\chi _{\mathrm{Ly}}$%
. So that $f_{\mathrm{Ly}}$ is continuous and hence bounded by compactness.

Hence boundedness is a necessary condition.

\vspace{12pt}%

\noindent%
\textbf{Example:} Let $\phi $ be a linear flow on a $d$-dimensional trivial
vector bundle $X\times V$ with two Lyapunov exponents $\lambda _{1}>\lambda
_{2}$ whose Oseledets subspaces have dimension $k$ and $d-k$. Then $H_{%
\mathrm{Ly}}$ is the diagonal matrix $\mathrm{diag}\{\lambda _{1},\ldots
,\lambda _{1},\lambda _{2},\ldots ,\lambda _{2}\}$ with $\lambda _{1}$
having multiplicity $k$. The Oseledets section is given by a map $f_{\mathrm{%
Ly}}:X\rightarrow \mathrm{Ad}\left( G\right) H_{\mathrm{Ly}}$ such that the
Oseledets subspaces at $x\in X$ are the eigenspaces $V_{\lambda _{1}}\left(
x\right) $ and $V_{\lambda _{2}}\left( x\right) $ of $f_{\mathrm{Ly}}\left(
x\right) $. To say that $f_{\mathrm{Ly}}$ is bounded means that the
subspaces $V_{\lambda _{1}}\left( x\right) $ and $V_{\lambda _{2}}\left(
x\right) $ have a positive distance.

\subsection{Refinement of the Lyapunov spectrum of other invariant measures
\label{cond2}}

Denote by $\mathcal{P}_{X}\left( \phi \right) $ the set of $\phi $-invariant
probability measures on $X$. For each ergodic measure $\rho \in \mathcal{P}%
_{X}\left( \phi \right) $ we have its Lyapunov spectrum $H_{\mathrm{Ly}%
}\left( \rho \right) $ and the corresponding flag type $\Theta _{\mathrm{Ly}%
}\left( \rho \right) =\{\alpha \in \Sigma :\alpha \left( H_{\mathrm{Ly}%
}\left( \rho \right) \right) =0\}$. Our second condition is

\begin{itemize}
\item $\Theta _{\mathrm{Ly}}\left( \rho \right) \subset \Theta _{\mathrm{Ly}%
}\left( \nu \right) $ for every ergodic $\rho \in \mathcal{P}_{X}\left( \phi
\right) $.
\end{itemize}

This is a necessary condition for $\Theta _{\mathrm{Ly}}\left( \nu \right)
=\Theta _{\mathrm{Mo}}$. To see this let $\rho \in \mathcal{P}_{X}\left(
\phi \right) $ be ergodic and denote by $Y\subset X$ its support. Let $%
\Theta _{\mathrm{Mo}}\left( Y\right) $ be the flag type of the Morse
decomposition of the flow restricted to $Y$ (that is, to the fibers above $Y$%
). We have $\Theta _{\mathrm{Mo}}\left( Y\right) \subset \Theta _{\mathrm{Mo}%
}\left( \nu \right) $ because the Morse components of the flow restricted to
$Y$ are contained in the components over $X$. However, $\Theta _{\mathrm{Ly}%
}\left( \rho \right) \subset \Theta _{\mathrm{Mo}}\left( Y\right) $, so that
the equality $\Theta _{\mathrm{Mo}}=\Theta _{\mathrm{Ly}}\left( \nu \right) $
implies
\[
\Theta _{\mathrm{Ly}}\left( \rho \right) \subset \Theta _{\mathrm{Mo}}\left(
Y\right) \subset \Theta _{\mathrm{Mo}}=\Theta _{\mathrm{Ly}}\left( \nu
\right) .
\]

\vspace{12pt}%

\noindent%
\textbf{Example:} Let $\phi $ be a linear flow on a $d$-dimensional vector
bundle $X\times V$ with Lyapunov spectrum $\lambda _{1}>\cdots >\lambda _{s}$
with multiplicities $k_{1},\ldots ,k_{s}$ with respect to $\nu $. Then this
condition means that the Lyapunov spectrum $\mu _{1}>\cdots >\mu _{t}$ with
respect to another ergodic measure $\rho $ have multiplicities $r_{1},\ldots
,r_{t}$ that satisfy $r_{1}+\cdots +r_{i_{1}}=k_{1}$, $r_{i_{1}}+\cdots
+r_{i_{2}}=k_{2}$, etc.

\subsection{Attracting and repeller measures\label{cond3}}

Recall that we defined an attractor measure on a partial flag manifold $%
\mathbb{E}_{\Theta }$ to be the projection of an ergodic invariant measure $%
\mu $ on $\mathbb{E}$ such that
\[
H_{\mathrm{Ly}}\left( \mu \right) =\int qd\mu \in \mathrm{cl}\mathfrak{a}%
^{+}.
\]

In the specific flag bundle $\mathbb{E}_{\Theta _{\mathrm{Ly}}}$ the
attractor measure $\mu _{\Theta _{\mathrm{Ly}}}^{+}$ is unique and has a
desintegration over $\nu $ by Dirac measures on the fibers above a set $%
\Omega $ of total measure $\nu $. We denote by $\mathrm{att}_{\Theta _{%
\mathrm{Ly}}}\left( \nu \right) $ the support of $\mu _{\Theta _{\mathrm{Ly}%
}}^{+}$.

Analogously in the dual flag bundle $\mathbb{E}_{\Theta _{\mathrm{Ly}}^{*}}$
there is a unique repeller measure $\mu _{\Theta _{\mathrm{Ly}}}^{-}$. We
denote by $\mathrm{rep}_{\Theta _{\mathrm{Ly}}^{*}}\left( \nu \right) $ the
support of $\mu _{\Theta _{\mathrm{Ly}}}^{-}$.

Let $\rho \in \mathcal{P}_{X}\left( \phi \right) $ be an ergodic measure
with support $Y\subset X$, which is an invariant subset. The subset $\pi
_{\Theta _{\mathrm{Ly}}}^{-1}\left( Y\right) \cap \mathrm{att}_{\Theta _{%
\mathrm{Ly}}}\left( \nu \right) $ is invariant as well. We denote by $%
\mathcal{E}_{\Theta _{\mathrm{Ly}}}^{+}\left( \rho \right) $ the set of
ergodic probability measures with support contained in $\pi _{\Theta _{%
\mathrm{Ly}}}^{-1}\left( Y\right) \cap \mathrm{att}_{\Theta _{\mathrm{Ly}%
}}\left( \nu \right) $ that project down to $\rho $. Also, we put $\mathcal{E%
}_{\Theta _{\mathrm{Ly}}^{*}}^{-}\left( \rho \right) $ for the set of
ergodic probability measures with support in $\pi _{\Theta _{\mathrm{Ly}%
}}^{-1}\left( Y\right) \cap \mathrm{rep}_{\Theta^{*}_{\mathrm{Ly}}}\left(
\nu \right) $ that project down to $\rho $. Both sets $\mathcal{E}_{\Theta _{%
\mathrm{Ly}}}^{+}\left( \rho \right) $ and $\mathcal{E}_{\Theta _{\mathrm{Ly}%
}^{*}}^{-}\left( \rho \right) $ are not empty.

Now we can state our third condition.

\begin{itemize}
\item Any $\theta \in \mathcal{E}_{\Theta _{\mathrm{Ly}}}^{+}\left( \rho
\right) $ is an attractor measure and any $\theta \in \mathcal{E}_{\Theta _{%
\mathrm{Ly}}^{*}}^{-}\left( \rho \right) $ is a reppeller measure for $\phi $%
.
\end{itemize}

If $\Theta _{\mathrm{Mo}}=\Theta _{\mathrm{Ly}}\left( \nu \right) $ then the
attractor Morse component $\mathcal{M}_{\Theta _{\mathrm{Ly}}}=\mathcal{M}%
_{\Theta _{\mathrm{Ly}}}\left( 1\right) $ in $\mathbb{E}_{\Theta _{\mathrm{Ly%
}}}$ is the image of a section $\xi :X\rightarrow \mathbb{E}_{\Theta _{%
\mathrm{Ly}}}$ and contains the support \textrm{att}$_{\Theta _{\mathrm{Ly}%
}}\left( \nu \right) $ of the attractor measure. This implies that $\mathcal{%
M}_{\Theta _{\mathrm{Ly}}}=\mathrm{att}_{\Theta _{\mathrm{Ly}}}\left( \nu
\right) $, so that $\mathcal{M}=\pi _{\Theta _{\mathrm{Ly}}}^{-1}\left(
\mathrm{att}_{\Theta _{\mathrm{Ly}}}\left( \nu \right) \right) $ is the
attractor Morse component $\mathcal{M}$ in the maximal flag bundle. Now the
Morse spectrum $\Lambda _{\mathrm{Mo}}\left( \mathcal{M}\right) $ of $%
\mathcal{M}$ is contained in the cone
\[
\mathfrak{a}_{\Theta _{\mathrm{Mo}}}^{+}=\{H\in \mathfrak{a}:\forall \alpha
\notin \langle \Theta _{\mathrm{Mo}}\rangle ,\,\alpha \left( H\right) >0\}
\]%
(see \cite{smsec}). Hence any Lyapunov exponent of $\mathcal{M}$ belongs to $%
\mathfrak{a}_{\Theta _{\mathrm{Mo}}}^{+}$. By projecting down to $\mathbb{E}%
_{\Theta _{\mathrm{Ly}}}$ the measures with support in $\mathcal{M}$ we see
that any $\theta $ with support in $\pi _{\Theta _{\mathrm{Ly}}}^{-1}\left(
Y\right) \cap \mathrm{att}_{\Theta _{\mathrm{Ly}}}\left( \nu \right) $ is an
attractor measure.

The same proof with the backward flow shows that  $\theta \in \mathcal{E}%
_{\Theta _{\mathrm{Ly}}^{\ast }}^{-}\left( \rho \right) $ is a reppeller
measure.

\subsection{Oseledets decompositions for other measures}

The second and third conditions above refer to ergodic measures $\rho $ on $%
X $ different from the initial measure $\nu $. These two conditions can be
summarized in just one condition on the Oseledets section for the ergodic
measures $\rho \in \mathcal{P}_{X}\left( \phi \right) $.

Given $\rho \in \mathcal{P}_{X}\left( \phi \right) $, write $\chi ^{\rho }$
for its Oseledets section and $\xi ^{\rho }$ and $\xi ^{\rho *}$ for the
corresponding sections on $\mathbb{E}_{\Theta _{\mathrm{Ly}}\left( \rho
\right) }$ and $\mathbb{E}_{\Theta _{\mathrm{Ly}}^{*}\left( \rho \right) }$,
respectively.

\begin{defi}
We say that $\chi ^{\rho }$ is contained in the Oseledets section of $\nu $
in case the two conditions are satisfied

\begin{enumerate}
\item $\Theta _{\mathrm{Ly}}\left( \rho \right) \subset \Theta _{\mathrm{Ly}%
}\left( \nu \right) $. In this case there is the fibration $p:Q\mathcal{O}%
_{\Theta _{\mathrm{Ly}}\left( \rho \right) }\rightarrow Q\mathcal{O}_{\Theta
_{\mathrm{Ly}}\left( \nu \right) }$.

\item $p\left( \mathrm{im}\chi ^{\rho }\right) \subset \mathrm{cl}\left(
\mathrm{im}\chi \right) $, where $\chi $ is the Oseledets section of $\nu $.
\end{enumerate}
\end{defi}

The second condition implies that the image of the sections $\xi ^{\rho }$
and $\xi ^{\rho *}$ project onto the $\mathrm{cl}\left( \mathrm{im}\xi
\right) $ and $\mathrm{cl}\left( \mathrm{im}\xi ^{*}\right) $, by the
fibrations $\mathbb{E}_{\Theta _{\mathrm{Ly}}\left( \rho \right)
}\rightarrow \mathbb{E}_{\Theta _{\mathrm{Ly}}\left( \nu \right) }$ and $%
\mathbb{E}_{\Theta _{\mathrm{Ly}}^{*}\left( \rho \right) }\rightarrow
\mathbb{E}_{\Theta _{\mathrm{Ly}}^{*}\left( \nu \right) }$, respectively.

Since the attractor and repeller measures for $\rho $ desintegrate according
to the sections $\xi ^{\rho }$ and $\xi ^{\rho *}$, respectively, it follows
that the second and third conditions above is equivalent to have $\chi
^{\rho }$ contained in $\chi $.

\section{Sufficience of the conditions}

\label{sufcon}

We apply here the main Lemma \ref{lemmaintech} to get sufficience of the
conditions of the last section and thus prove the following characterization
for the equality of Morse and Oseledets decompositions.

\begin{teorema}
\label{teomainequal}Suppose the invariant measure on the base space is
ergodic. Then the three conditions together --- bounded section (\ref{cond1}%
), refinement of Lyapunov spectrum (\ref{cond2}) and attractor-repeller
measures (\ref{cond3}) --- are necessary and sufficient to have $\Theta _{%
\mathrm{Ly}}=\Theta _{\mathrm{Mo}}$ and $\chi _{\mathrm{Ly}}=\chi _{\mathrm{%
Mo}}$.
\end{teorema}

As before we have the sections $\xi :\Omega \rightarrow \mathbb{E}_{\Theta _{%
\mathrm{Ly}}}$ and $\xi ^{*}:\Omega \rightarrow \mathbb{E}_{\Theta _{\mathrm{%
Ly}}^{*}}$, respectively, that are combined to give the Oseledets section $%
\chi _{\mathrm{Ly}}:\Omega \subset X\rightarrow \mathcal{A}_{\Theta _{%
\mathrm{Ly}}}$.

We apply Lemma \ref{lemmaintech} with

\begin{enumerate}
\item $A=\mathrm{cl}\left( \mathrm{im}\xi \right) $, which is the support of
the unique attrator measure $\mu _{\Theta _{\mathrm{Ly}}}^{+}$ in $\mathbb{E}%
_{\Theta _{\mathrm{Ly}}}$.

\item $C=\mathrm{cl}\left( \mathrm{im}\xi ^{*}\right) $, which is the
support of the unique repeller measure $\mu _{\Theta _{\mathrm{Ly}}}^{-}$ in
$\mathbb{E}_{\Theta _{\mathrm{Ly}}^{*}}$, and

\item $B=\mathrm{cl}\bigcup_{w\neq 1}\mathrm{st}_{\Theta _{\mathrm{Ly}%
}}\left( x,w\right) $, $x\in \Omega $. That is, $B$ is the closure of the
set of elements that are \textbf{not} transversal to $\xi ^{*}\left(
x\right) $, $x\in \Omega $.
\end{enumerate}

Alternatively we have the following characterization of $B$ in terms of the
closure of the dual section $\xi ^{*}$.

\begin{proposicao}
\label{propbtransv}An element $v\in \mathbb{E}_{\Theta _{\mathrm{Ly}}}$
belongs to $B$ if and only if it is not transversal to some $w\in \mathrm{cl}%
\left( \mathrm{im}\xi ^{*}\right) $ in the same fiber as $v$.
\end{proposicao}

\begin{profe}
Take  local trivializations so that locally $\mathbb{E}_{\Theta _{\mathrm{Ly%
}}}\simeq U\times \mathbb{F}_{\Theta _{\mathrm{Ly}}}$, $\mathbb{E}_{\Theta _{%
\mathrm{Ly}}^{*}} \simeq U\times \mathbb{F}_{\Theta _{\mathrm{Ly}}^{*}}$ ($%
U\subset X$ open), $\xi :U\rightarrow \mathbb{F}_{\Theta _{\mathrm{Ly}}}$
and $\xi ^{*}:U\rightarrow \mathbb{F}_{\Theta _{\mathrm{Ly}}^{*}}$. If $%
v=\left( x,b\right) \in B$ then there exists a sequence $\left(
x_{n},b_{n}\right) \rightarrow v$ with $b_{n}$ not transversal to $\xi
^{*}\left( x_{n}\right) $. By taking a subsequence we can assume that $\xi
^{*}\left( x_{n}\right) $ converges to $b^{*}\in \mathbb{F}_{\Theta _{%
\mathrm{Ly}}^{*}}$. Then the pair $\left( b_{n},\xi ^{*}\left( x_{n}\right)
\right) $ converges to $\left( b,b^{*}\right) \in \mathbb{F}_{\Theta _{%
\mathrm{Ly}}}\times \mathbb{F}_{\Theta _{\mathrm{Ly}}^{*}}$. Now, the set of
nontransversal pairs in $\mathbb{F}_{\Theta _{\mathrm{Ly}}}\times \mathbb{F}%
_{\Theta _{\mathrm{Ly}}^{*}}$ is closed. Hence $b$ and $b^{*}$ are not
transversal, showing that $v=\left( x,b\right) $ is not transversal to $%
w=\left( x,b^{*}\right) \in \mathrm{cl}\left( \mathrm{im}\xi ^{*}\right) $.

Conversely, suppose that $v=\left( x,b\right) \in \mathbb{E}_{\Theta _{%
\mathrm{Ly}}}$ is not tranversal to $w=\left( x,b^{*}\right) \in \mathrm{cl}%
\left( \mathrm{im}\xi ^{*}\right) $. Then $b^{*}=\lim \xi ^{*}\left(
x_{n}\right) $ with $\lim x_{n}=x$. By Lemma \ref{lemseqtransv} there exists
a sequence $b_{n}\in \mathbb{F}_{\Theta _{\mathrm{Ly}}}$ such that $b_{n}$
is not transversal to $\xi ^{*}\left( x_{n}\right) $ and $\lim b_{n}=b$.
Hence $\left( x_{n},b_{n}\right) \in B$ and $\lim \left( x_{n},b_{n}\right)
=\left( x,b\right) =v$, showing that $v\in B$.
\end{profe}

Clearly, $A$, $B$ and $C$ are compact sets. Also, $A$ and $C$ are invariant
because the sections $\xi $ and $\xi ^{*}$ are invariant, and since
transversality is preserved by the flow, it follows that $B$ is invariant as
well.

Now we verify that the assumptions of Lemma \ref{lemmaintech} hold in
presence of the three conditions of Theorem \ref{teomainequal}. Statements
(3) and (4) of Lemma \ref{lemmaintech} are the same as the refinement of
Lyapunov spectrum and attractor-repeller measures conditions, respectively.
Item (2) of Lemma \ref{lemmaintech} is the above proposition. So it remains
to prove that $A$ and $B$ are disjoint. This is the only place where the
bounded condition is used.

\begin{proposicao}
$A\cap B=\emptyset $. Precisely, if $v\in A$ and $w\in \mathrm{cl}\left(
\mathrm{im}\xi ^{*}\right) $ then $v$ and $w$ are transversal, and if $v\in
B $ then there exists $w\in \mathrm{cl}\left( \mathrm{im}\xi ^{*}\right) $
in the same fiber which is not transversal to $v$.
\end{proposicao}

\begin{profe}
Since the restriction of $f_{\mathrm{Ly}}$ to $R_{\Omega }$ is bounded its
image in $\mathcal{O}_{\mathrm{Ly}}=\mathrm{Ad}\left( G\right) H_{\mathrm{Ly}%
}$ has compact closure. By the equality $\chi _{\mathrm{Ly}}\left( x\right)
=p\cdot f_{\mathrm{Ly}}\left( p\right) $ ($p\in Q$ with $\pi (p)=x$) it
follows that $\mathrm{cl}\left( \mathrm{im}\chi _{\mathrm{Ly}}\right) $ is a
compact subset of the bundle $Q\times _{G}\mathcal{O}_{\mathrm{Ly}}$. After
identifying $\mathcal{O}_{\mathrm{Ly}}$ with an open subset of $\mathbb{F}%
_{\Theta _{\mathrm{Ly}}}\times \mathbb{F}_{\Theta _{\mathrm{Ly}}^{\ast }}$
we get a section of $Q\times _{G}\left( \mathbb{F}_{\Theta _{\mathrm{Ly}%
}}\times \mathbb{F}_{\Theta _{\mathrm{Ly}}^{\ast }}\right) $ over $\Omega $
also denoted by $\chi _{\mathrm{Ly}}$. The image of this section is
contained in the open subset of those pairs in $Q\times _{G}\left( \mathbb{F}%
_{\Theta _{\mathrm{Ly}}}\times \mathbb{F}_{\Theta _{\mathrm{Ly}}^{\ast
}}\right) $ that are transversal to each other. By compactness the closure
of the image of $\chi _{\mathrm{Ly}}$ contains also only transversal pairs.

Now let $p:Q\times _{G}\left( \mathbb{F}_{\Theta _{\mathrm{Ly}}}\times
\mathbb{F}_{\Theta _{\mathrm{Ly}}^{*}}\right) \rightarrow \mathbb{E}_{\Theta
_{\mathrm{Ly}}}$ and $p^{*}:Q\times _{G}\left( \mathbb{F}_{\Theta _{\mathrm{%
Ly}}}\times \mathbb{F}_{\Theta _{\mathrm{Ly}}^{*}}\right) \rightarrow
\mathbb{E}_{\Theta _{\mathrm{Ly}}^{*}}$ be the canonical projections. Then
\[
\xi =p\circ \chi _{\mathrm{Ly}}\qquad \mathrm{and}\qquad \xi ^{*}=p^{*}\circ
\chi _{\mathrm{Ly}}.
\]
Hence by compactness $p\left( \mathrm{cl}\left( \mathrm{im}\chi _{\mathrm{Ly}%
}\right) \right) =\mathrm{cl}\left( \mathrm{im}\xi \right) $ and $%
p^{*}\left( \mathrm{cl}\left( \mathrm{im}\chi _{\mathrm{Ly}}\right) \right) =%
\mathrm{cl}\left( \mathrm{im}\xi ^{*}\right) $. It follows that two elements
$v\in A=\mathrm{cl}\left( \mathrm{im}\xi \right) $ and $w\in \mathrm{cl}%
\left( \mathrm{im}\xi ^{*}\right) $ are transversal to each other, if they
are in the same fiber.

On the other hand if $v\in B$ then by Proposition \ref{propbtransv}, there
exists $w\in \mathrm{cl}\left( \mathrm{im}\xi ^{*}\right) $ such that $v$
and $w$ are in the same fiber and are not transversal. Hence $v\notin A$,
concluding that $A$ and $B$ are disjoint.
\end{profe}

\vspace{12pt}%

\noindent \textbf{End of proof of Theorem \ref{teomainequal}}: By Lemma \ref%
{lemmaintech}, $A$ and $B$ define a Morse decomposition of $\mathbb{E}%
_{\Theta _{\mathrm{Ly}}}$ with $A$ the attractor component. Hence $A=\mathrm{%
cl}\left( \mathrm{im}\xi \right) $ contains the unique attractor component $%
\mathcal{M}_{\Theta _{\mathrm{Ly}}}^{+}$ of the finest Morse decomposition
of $\mathbb{E}_{\Theta _{\mathrm{Ly}}}$. On the other hand by Proposition %
\ref{prooselcontmor} the Oseledets component $\mathrm{im}\xi \subset
\mathcal{M}_{\Theta _{\mathrm{Ly}}}^{+}$. Therefore $A=\mathrm{cl}\left(
\mathrm{im}\xi \right) \subset \mathcal{M}_{\Theta _{\mathrm{Ly}}}^{+}$, so
that they are equal.

\section{Uniquely ergodic base spaces\label{secunique}}

When the flow on the base space has unique invariant (and hence ergodic)
probability measure $\nu $ the second and third conditions of Section \ref%
{threecon} are meaningless. Hence, in this case, a necessary and sufficient
condition to have equality of Oseledets and Morse decompositions is that the
Oseledets section for $\nu $ is bounded (first condition of Section \ref%
{threecon}).

From another point of view the Morse spectrum $\Lambda _{\mathrm{Mo}}$ of
the attractor component $\mathcal{M}^{+}$ is a compact convex set whose
extremal points are Lyapunov exponents given by integrals with respect to
invariant measures on the maximal flag bundle. By the results of Section \ref%
{secmedinv} any such integral Lyapunov exponent is a regular Lyapunov
exponent of an invariant probability in the base space. Just one of these
Lyapunov exponents belongs to \textrm{cl}$\mathfrak{a}^{+}$, which is the
polar exponent $H_{\mathrm{Ly}}$ associated to the measure.

Hence, $\Lambda _{\mathrm{Mo}}$ has a unique extremal point in \textrm{cl}$%
\mathfrak{a}^{+}$ if the flow on the base space is uniquely ergodic.

\begin{proposicao}
Suppose that the flow on the base space $X$ has a unique invariant
probability measure $\nu $ with $\mathrm{supp}\nu =X$. Let $H_{\mathrm{Ly}%
}=H_{\mathrm{Ly}}\left( \nu \right) $ be its polar exponent. Then $\Lambda _{%
\mathrm{Mo}}$ is the polyhedron whose vertices are $wH_{\mathrm{Ly}}$, $w\in
\mathcal{W}_{\Theta _{\mathrm{Mo}}}$.
\end{proposicao}

\begin{profe}
Since the convex set $\Lambda _{\mathrm{Mo}}$ is invariant by $\mathcal{W}%
_{\Theta _{\mathrm{Mo}}}$ and $H_{\mathrm{Ly}}\in \Lambda _{\mathrm{Mo}}$ we
have that the polyhedron with vertices in $\mathcal{W}_{\Theta _{\mathrm{Mo}%
}}\left( H_{\mathrm{Ly}}\right) $ is contained in $\Lambda _{\mathrm{Mo}}$.
Conversely, suppose that $H$ is an extremal point of $\Lambda _{\mathrm{Mo}}$%
. Then there exists $w\in \mathcal{W}$ such that $wH\in \mathrm{cl}\mathfrak{%
a}^{+}$. We claim that $w\in \mathcal{W}_{\Theta _{\mathrm{Mo}}}$. In fact,
by Weyl group invariance of the Morse spectrum there exists a Morse
component $\mathcal{M}$ such that $wH\in \Lambda _{\mathrm{Mo}}\left(
\mathcal{M}\right) $ (see \cite{smsec}). But the attractor component $%
\mathcal{M}^{+}$ is the only one whose Morse spectrum meets $\mathrm{cl}%
\mathfrak{a}^{+}$. So that $\mathcal{M=M}^{+}$ and $wH\in \Lambda _{\mathrm{%
Mo}}=\Lambda _{\mathrm{Mo}}\left( \mathcal{M}^{+}\right) $ and since the
spectra of distinct Morse components are disjoint we have $w\Lambda _{%
\mathrm{Mo}}=\Lambda _{\mathrm{Mo}}$, implying that $w\in \mathcal{W}%
_{\Theta _{\mathrm{Mo}}}$.

Now $wH\in \mathrm{cl}\mathfrak{a}^{+}$ is an extremal point of $\Lambda _{%
\mathrm{Mo}}=w\Lambda _{\mathrm{Mo}}$ hence $wH=H_{\mathrm{Ly}}$. Therefore,
$\mathcal{W}_{\Theta _{\mathrm{Mo}}}H_{\mathrm{Ly}}$ is the set of extremal
points of $\Lambda _{\mathrm{Mo}}$, concluding the proof.
\end{profe}

\begin{teorema}
Suppose that the flow on the base space $X$ has a unique invariant
probability measure $\nu $ with $\mathrm{supp}\nu =X$. Then the following
conditions are equivalent.

\begin{enumerate}
\item $\Theta _{\mathrm{Ly}}=\Theta _{\mathrm{Mo}}$.

\item The Oseledets section for $\nu $ is bounded.

\item $\alpha \left( \Lambda _{\mathrm{Mo}}\right) =\{0\}$ for all $\alpha
\in \Theta _{\mathrm{Mo}}$.
\end{enumerate}

If these conditions hold then $\Lambda _{\mathrm{Mo}}=\{H_{\mathrm{Ly}}\}$.
\end{teorema}

\begin{profe}
As mentioned above the equivalence between the first two conditions is a
consequence of the main Theorem \ref{teomainequal} and the fact that $\nu $
is the only ergodic measure on $X$. Now, $\Theta _{\mathrm{Ly}}=\Theta _{%
\mathrm{Mo}}$ means that $\alpha \left( H_{\mathrm{Ly}}\right) =0$ for all $%
\alpha \in \Theta _{\mathrm{Mo}}$. If this happens then any $\alpha \in
\Theta _{\mathrm{Mo}}$ annihilates on the polyhedron with vertices $wH_{%
\mathrm{Ly}}$, $w\in \mathcal{W}_{\Theta _{\mathrm{Mo}}}$. Hence $\alpha \in
\Theta _{\mathrm{Mo}}$ is zero on $\Lambda _{\mathrm{Mo}}$ by the above
proposition. Conversely, if (3) holds then $\alpha \left( H_{\mathrm{Ly}%
}\right) =0$ for all $\alpha \in \Theta _{\mathrm{Mo}}$, because $H_{\mathrm{%
Ly}}\in \Lambda _{\mathrm{Mo}}$.

Finally, if $\alpha \left( H_{\mathrm{Ly}}\right) =0$ for all $\alpha \in
\Theta _{\mathrm{Mo}}$ then $wH_{\mathrm{Ly}}=H_{\mathrm{Ly}}$ for every $%
w\in \mathcal{W}_{\Theta _{\mathrm{Mo}}}$ so that $\Lambda _{\mathrm{Mo}%
}=\{H_{\mathrm{Ly}}\}$ by the previous proposition.
\end{profe}

By piecing together known results in the literature we can have examples of
flows over uniquely ergodic systems for which $\Theta _{\mathrm{Ly}}\neq
\Theta _{\mathrm{Mo}}$. Indeed as proved by Furman \cite{fur} the Lyapunov
spectrum is discontinuous at a non-uniform cocycle with values in $\mathrm{Gl%
}\left( d,\mathbb{R}\right) $, that is, at a flow on the trivial bundle $%
X\times \mathrm{Gl}\left( d,\mathbb{R}\right) $. The result of \cite{fur}
(see Theorem 5) makes the assumption that the flow on the base space is
equicontinuous, which is satisfied, e.g., by the translations on compact
groups, like an irrational rotation on the circle $S^{1}$.

Now in Herman \cite{her} there is an example of a non-uniform cocycle with
values in $\mathrm{Sl}\left( 2,\mathbb{R}\right) $ over the irrational
rotation. Thus that example is a discontinuity point of the Lyapunov
spectrum. Finally, in Section \ref{seccont} below we prove continuity of the
whole Lyapunov spectrum if $\Theta _{\mathrm{Ly}}=\Theta _{\mathrm{Mo}}$.
Hence, we get $\Theta _{\mathrm{Ly}}\neq \Theta _{\mathrm{Mo}}$ for the
example in \cite{her}.

\section{Product of i.i.d. sequences\label{seciid}}

The product of independent identically distributed (i.i.d.) random elements
in $G$ yield flows on product spaces $\mathcal{X}\times G$ with plenty of
invariant measures on $\mathcal{X}$. In this section we provide an example
such flow that violates the second condition of Section \ref{threecon} and
hence has distinct Morse and Oseledets decompositions.

Let $C\subset G$ be a compact subset and form the product $\mathcal{X}=C^{%
\mathbb{Z}}$ endowed with the compact product topology. The shift $\tau
\left( \left( x_{n}\right) \right) =\left( x_{n+1}\right) _{n\in \mathbb{Z}}$
is a homeomorphism and hence defines a continuous flow on $\mathcal{X}$.

Now, let $\mu $ be a probability measure with $\mathrm{supp}\mu =C$ and take
the product measure $\mu ^{\times \mathbb{Z}}$ on $\mathcal{X}$. Then $\mu
^{\times \mathbb{Z}}$ is ergodic with respect to the shift $\tau $ and $%
\mathrm{supp}\mu ^{\times \mathbb{Z}}=\mathcal{X}$.

These data defines the continuous flow $\phi _{n}^{\mu }$ on $\mathcal{X}%
\times G$ by $\phi _{n}^{\mu }\left( \mathbf{x},g\right) =\left( \tau
^{n}\left( \mathbf{x}\right) ,\rho ^{\mu }\left( n,\mathbf{x}\right)
g\right) $, where $\mathbf{x}=\left( x_{n}\right) _{n\in \mathbb{Z}}\in C^{%
\mathbb{Z}}\subset G^{\mathbb{Z}}$ and
\[
\rho ^{\mu }\left( n,\mathbf{x}\right) =\left\{
\begin{array}{lll}
x_{n-1}\cdots x_{0} & \mathrm{se} & n\geq 0 \\
x_{1}^{-1}\cdots x_{n}^{-1} & \mathrm{se} & n<0.%
\end{array}%
\right.
\]

The $\mathfrak{a}$-Lyapunov spectrum of $\phi ^{\mu }$ were finded by
Guivarch'-Raugi \cite{gr}. To state their result we recall the following
concepts.

\begin{enumerate}
\item A subgroup $H\subset G$ is \textit{totally irreducible} if it does not
leave invariant a subset which is a finite union of complements of Bruhat
cells in their respective closures (Schubert cells).

\item A sequence $g_{n}\in G$ is said to be \textit{contracting} with
respect to the maximal flag manifold if its polar decomposition $%
g_{n}=u_{n}h_{n}v_{n}\in K\left( \mathrm{cl}A^{+}\right) K$ is such that
\[
\lim_{n\rightarrow \infty }\alpha \left( \log h_{n}\right) =\infty
\]%
for every positive root $\alpha $.
\end{enumerate}

Now denote by $G_{\mu }$ and $S_{\mu }$ the subgroup and semigroup generated
by $\mathrm{supp}\mu =C$, respectively. Then we have the following result of
\cite{gr}, Theorem 2.6.

\begin{teorema}
Let $\mu $ be a probability measure on $G$, and suppose that

\begin{enumerate}
\item the subgroup $G_{\mu }$ is totally irreducible.

\item The semigroup $S_{\mu }$ has a contracting sequence with respect to $%
\mathbb{F}$.
\end{enumerate}

Then the polar exponent of $\phi ^{\mu }$ is regular, that is, belongs to $%
\mathfrak{a}^{+}$. This means that $\Theta _{\mathrm{Ly}}\left( \mu ^{\times
\mathbb{Z}}\right) =\emptyset $.
\end{teorema}

Both conditions of this theorem are satisfied if $S_{\mu }$ has nonempty
interior in $G$:

\begin{enumerate}
\item If $\mathrm{int}S_{\mu }\neq \emptyset $ then $G_{\mu }=G$ because $G$
is assumed to be connected.

\item If $\mathrm{int}S_{\mu }\neq \emptyset $ then there exists a regular $%
h\in \mathrm{int}S_{\mu }$ (see \cite{smics}, Lemma 3.2). Then $h^{n}\in
S_{\mu }$ is a contracting sequence with respect to $\mathbb{F}$.
\end{enumerate}

Hence we get the following consequence of Guivarch'-Raugi \cite{gr} result.

\begin{corolario}
If $\mathrm{int}S_{\mu }\neq \emptyset $ then $\Theta _{\mathrm{Ly}}\left(
\phi ^{\mu }\right) =\emptyset $.
\end{corolario}

Now it is easy to give an exemple that does not satisfy the second condition
of Section \ref{threecon} and hence $\Theta _{\mathrm{Ly}}\left( \phi ^{\mu
}\right) \neq \Theta _{\mathrm{Mo}}\left( \phi ^{\mu }\right) $. In fact,
take a a nonregular element $h=\exp H\in \mathrm{cl}A^{+}$, that is, $\Theta
\left( H\right) =\{\alpha \in \Sigma :\alpha \left( H\right) =0\}\neq
\emptyset $ and a probability $\mu $ whose support $C$ contains $h$ in its
interior. For instance
\[
\mu =\frac{1}{I}f\cdot \eta
\]%
where $\eta $ is Haar measure and $I=\int_{G}f\left( g\right) \eta \left(
dg\right) <\infty $ with $f:G\rightarrow \mathbb{R}$ a nonnegative function
with $\mathrm{supp}f=C$.

Let $\rho =\delta _{\mathbf{x}_{h}}$ be the Dirac measure at\ the constant
sequence $\mathbf{x}_{h}=\left( x_{n}\right) _{n\in \mathbb{Z}}$, $x_{n}=h$.
Clearly $\rho $ is $\tau $-invariant and ergodic. Since
\[
\phi _{n}^{\mu }\left( \mathbf{x}_{h},1\right) =\left( \tau ^{n}\left(
\mathbf{x}_{h}\right) ,\rho \left( n,\mathbf{x}_{h}\right) \right) =\left(
\mathbf{x}_{h},h^{n}\right)
\]%
the polar exponent of $\delta _{\mathbf{x}_{h}}$ is
\[
\lim_{n\rightarrow +\infty }\frac{1}{n}\log h^{n}=\log h=H
\]%
which is not regular. Hence $\Theta _{\mathrm{Ly}}\left( \delta _{\mathbf{x}%
_{h}}\right) =\Theta \left( H\right) \neq \emptyset $ is not contained in $%
\Theta _{\mathrm{Ly}}\left( \mu ^{\times \mathbb{Z}}\right) =\emptyset $.
Therefore, the second condition of Section \ref{threecon} is violated and
the flag types $\Theta _{\mathrm{Ly}}\left( \mu ^{\times \mu }\right) $ and $%
\Theta _{\mathrm{Mo}}\left( \phi ^{\mu }\right) $ are different.

\section{Continuity of the Lyapunov spectrum\label{seccont}}

In this section we apply the differentiability result of \cite{ferrsm} to
show that the equality $\Theta _{\mathrm{Ly}}=\Theta _{\mathrm{Mo}}$ implies
continuity of the Lyapunov spectrum by perturbations of the original $\phi $
that do not change the flow on the base space.

Let $\mathcal{G}=\mathcal{G}\left( Q\right) $ be the gauge group of $%
Q\rightarrow X$, that is, the group of automorphisms of $Q$ that project to
the identity map of $X$. It is well known that $\mathcal{G}$ is a Banach Lie
group.

If $\sigma \in \mathcal{G}$ then $\phi $ and $\sigma \phi $ induce the same
map on $X$ and hence have the same ergodic measure $\nu $. Denote by $H_{%
\mathrm{Ly}}^{\sigma \phi }$ the polar spectrum $\sigma \phi $ with respect
to $\nu $. Assume as before that $\nu $ has full support. Then we have the
following continuity result.

\begin{teorema}
\label{teocont}If $\Theta _{\mathrm{Ly}}\left( \phi \right) =\Theta _{%
\mathrm{Mo}}$ then the map $\sigma \in \mathcal{G}\mapsto H_{\mathrm{Ly}%
}^{\sigma \phi }\in \mathrm{cl}\mathfrak{a}^{+}$ is continuous at $\sigma =%
\mathrm{id}$.
\end{teorema}

We work out separetely the proof for $\mathrm{Sl}\left( n,\mathbb{R}\right) $
in order to explain it in concrete terms. For this group $\mathfrak{a}$ is
the algebra of zero trace diagonal matrices and $\mathfrak{a}^{+}$ are those
with strictly decreasing eigenvalues. The simple set of roots is $\Sigma
=\{\alpha _{1},\ldots ,\alpha _{n-1}\}$ where $\alpha _{i}=\alpha
_{i,i+1}=\lambda _{i}-\lambda _{i+1}$ and $\lambda _{i}\in \mathfrak{a}%
^{\ast }$ maps the diagonal $H\in \mathfrak{a}$ to its $i$-th diagonal
entry. We denote by $\Delta =\{\delta _{1},\ldots ,\delta _{n-1}\}$ the set
of fundamental weights, which is defined by
\[
\frac{2\langle \alpha _{i},\delta _{j}\rangle }{\langle \alpha _{i},\alpha
_{i}\rangle }=\delta _{ij}
\]%
and is given by $\delta _{j}=\lambda _{1}+\cdots +\lambda _{j}$. The Morse
decomposition of $\phi $ on the flag bundles are determined by the subset $%
\Theta _{\mathrm{Mo}}\subset \Sigma $. Alternatively we can look at the
partition
\[
\{1,\ldots ,n\}=\{1,\ldots ,r_{1}\}\cup \{r_{1}+1,\ldots ,r_{2}\}\cup \cdots
\cup \{r_{k}+1,\ldots ,n\}
\]%
where $\Sigma \setminus $ $\Theta _{\mathrm{Mo}}=\{\alpha _{r_{1}},\ldots
,\alpha _{r_{k}}\}$. From the partition we recover $\Theta _{\mathrm{Mo}}$
as the set of $\alpha _{j,j+1}$ such that if $\left[ r,s\right] $ is the
interval of the partition containing $j$ then $r\leq j<s$.

For $H\in \mathrm{cl}\mathfrak{a}^{+}$ with $\alpha \left( H\right) =0$ for
all $\alpha \in \Theta _{\mathrm{Mo}}$ its eigenvalues $a_{i}$ are such that
$a_{i}=a_{j}$ if the indices $i,j$ belong to the same set of the partition.
If furthermore $H$ is such that $\Theta _{\mathrm{Mo}}=\{\alpha \in \Sigma
:\alpha \left( H\right) =0\}$ then the multiplicities of the eigenvalues of $%
H$ are the sizes of the sets of the partition.

Hence $\Theta _{\mathrm{Ly}}=\Theta _{\mathrm{Mo}}$ means that the
multiplicities of the Lyapunov exponents are given by the partition
associated to $\Theta _{\mathrm{Mo}}$.

Now, it was proved in \cite{ferrsm} that the map
\[
\sigma \in \mathcal{G}\mapsto \delta _{j}\left( H_{\mathrm{Ly}}^{\sigma \phi
}\right) \in \mathbb{R}_{+}
\]%
is differentiable at $\sigma =\mathrm{id}$ for any index $j$ such that $%
\alpha _{j,j+1}\notin \Theta _{\mathrm{Mo}}$. Since $\Delta $ is a basis of $%
\mathfrak{a}^{\ast }$ we get continuity of $H_{\mathrm{Ly}}^{\sigma \phi }$
if we prove that $\delta _{j}\left( H_{\mathrm{Ly}}^{\sigma \phi }\right) $
is continuous when $\alpha _{j,j+1}\in \Theta _{\mathrm{Mo}}$.

For this purpose we recall from \cite{alvsm} that $\delta _{j}\left( H_{%
\mathrm{Ly}}^{\sigma \phi }\right) $ is obtained as a limit furnished by the
sub-additive ergodic theorem. Namely,
\begin{equation}
\delta _{j}\left( H_{\mathrm{Ly}}^{\sigma \phi }\right) =\lim \frac{1}{k}%
\delta _{j}\left( \mathsf{a}_{\sigma }^{+}\left( k,x\right) \right)
=\inf_{k\geq 1}\frac{1}{k}\int \delta _{j}\left( \mathsf{a}_{\sigma
}^{+}\left( k,x\right) \right) \nu \left( dx\right)  \label{forlimsubad}
\end{equation}%
where $\mathsf{a}_{\sigma }^{+}\left( k,x\right) $, $x\in X$, is the polar
component of the flow define by $\sigma \phi $. (See Section 3.2 in \cite%
{alvsm}. Since $\delta _{j}$ is a fundamental weight $\delta _{j}\left(
\mathsf{a}_{\sigma }^{+}\left( k,x\right) \right) $ is a sub-additive
cocycle on the base space. As showed in \cite{alvsm} this cocycle can be
written as a norm in the space of a representation of $G$, which in this
case is the $j$-fold exterior product of $\mathbb{R}^{n}$.)

By (\ref{forlimsubad}) we have that $\sigma \mapsto \delta _{j}\left( H_{%
\mathrm{Ly}}^{\sigma \phi }\right) $ is upper semi-continuous.

To prove continuity take $j$ with $\alpha _{j,j+1}\in \Theta _{\mathrm{Mo}}$
and let $\left[ r,s\right] $, $r\leq j<s$, be the interval of the partition
that contains $j$. Assume by contradiction that there exist $c>0$ and a
sequence $\sigma _{k}\in \mathcal{G}$ converging to $\mathrm{id}$ such that
\[
\delta _{j}\left( H_{\mathrm{Ly}}^{\sigma _{k}\phi }\right) <\delta
_{j}\left( H_{\mathrm{Ly}}^{\phi }\right) -c.
\]%
Then we have two cases:

\begin{enumerate}
\item $s<n$. Then $\alpha _{s,s+1}\notin \Theta _{\mathrm{Mo}}$, so that $%
\sigma \mapsto \delta _{s}\left( H_{\mathrm{Ly}}^{\sigma \phi }\right) $ is
continuous. The same way $\delta _{r-1}\left( H_{\mathrm{Ly}}^{\sigma \phi
}\right) $ is continuous (where $\delta _{r-1}=0$ if $r=1$). Then for large $%
k$ we have
\[
\delta _{r-1}\left( H_{\mathrm{Ly}}^{\sigma _{k}\phi }\right) >\delta
_{r-1}\left( H_{\mathrm{Ly}}^{\phi }\right) -c/2.
\]%
Since $\lambda _{i_{1}}\geq \lambda _{i_{2}}$ on $\mathrm{cl}\mathfrak{a}%
^{+} $ if $i_{1}\leq i_{2}$ and the polar exponents $H_{\mathrm{Ly}}^{\sigma
_{k}\phi }\in \mathrm{cl}\mathfrak{a}^{+}$ we get
\[
\delta _{j}\left( H_{\mathrm{Ly}}^{\phi }\right) -c>\delta _{r-1}\left( H_{%
\mathrm{Ly}}^{\sigma _{k}\phi }\right) +\left( j-r+1\right) \lambda
_{j}\left( H_{\mathrm{Ly}}^{\sigma _{k}\phi }\right) .
\]%
Hence for large $k$ we have
\[
\delta _{j}\left( H_{\mathrm{Ly}}^{\phi }\right) -c>\delta _{r-1}\left( H_{%
\mathrm{Ly}}^{\phi }\right) -c/2+\left( j-r+1\right) \lambda _{j}\left( H_{%
\mathrm{Ly}}^{\sigma _{k}\phi }\right)
\]%
that is,
\[
\lambda _{j}\left( H_{\mathrm{Ly}}^{\sigma _{k}\phi }\right) <\frac{1}{%
\left( j-r+1\right) }\left( \delta _{j}\left( H_{\mathrm{Ly}}^{\phi }\right)
-\delta _{r-1}\left( H_{\mathrm{Ly}}^{\phi }\right) -c/2\right) .
\]

By the inequality $\delta _{s}=\delta _{j}+\lambda _{j+1}+\cdots +\lambda
_{s}\leq \delta _{j}+\left( s-j\right) \lambda _{j}$ that holds on $\mathrm{%
cl}\mathfrak{a}^{+}$ we get
\begin{equation}
\delta _{s}\left( H_{\mathrm{Ly}}^{\sigma _{k}\phi }\right) \leq \delta
_{j}\left( H_{\mathrm{Ly}}^{\phi }\right) -c+\frac{s-j}{j-r+1}\left( \delta
_{j}\left( H_{\mathrm{Ly}}^{\phi }\right) -\delta _{r-1}\left( H_{\mathrm{Ly}%
}^{\phi }\right) -c/2\right) .  \label{fordesigualdade}
\end{equation}%
Now we use the assumption $\Theta _{\mathrm{Ly}}\left( \phi \right) =\Theta
_{\mathrm{Mo}}$ which implies that $\delta _{j}\left( H_{\mathrm{Ly}}^{\phi
}\right) =\delta _{r-1}\left( H_{\mathrm{Ly}}^{\phi }\right) +\left(
j-r+1\right) \lambda _{j}\left( H_{\mathrm{Ly}}^{\phi }\right) $ and $\delta
_{s}\left( H_{\mathrm{Ly}}^{\phi }\right) =\delta _{j}\left( H_{\mathrm{Ly}%
}\left( \phi \right) \right) +\left( s-j\right) \lambda _{j}\left( H_{%
\mathrm{Ly}}^{\phi }\right) $. Hence the last term in (\ref{fordesigualdade}%
) becomes
\[
\frac{s-j}{j-r+1}\left( \left( j-r+1\right) \lambda _{j}\left( H_{\mathrm{Ly}%
}^{\phi }\right) -c/2\right) =\left( s-j\right) \lambda _{j}\left( H_{%
\mathrm{Ly}}^{\phi }\right) -\frac{s-j}{j-r+1}\frac{c}{2},
\]%
so that for large $k$ it holds
\begin{equation}
\delta _{s}\left( H_{\mathrm{Ly}}^{\sigma _{k}\phi }\right) \leq \delta
_{s}\left( H_{\mathrm{Ly}}^{\phi }\right) -\frac{s-j}{j-r+1}\frac{c}{2}
\label{fordesigualfinal}
\end{equation}%
which contradicts the continuity of $\delta _{s}\left( H_{\mathrm{Ly}%
}^{\sigma \phi }\right) $.

\item $s=n$. If $r=1$ then $\Theta _{\mathrm{Ly}}\left( \phi \right) =\Theta
_{\mathrm{Mo}}=\Sigma $ so that $H_{\mathrm{Ly}}^{\phi }=0$ and continuity
follows by upper semi-continuity. When $r\neq 1$ we get continuity of $%
\delta _{r-1}\left( H_{\mathrm{Ly}}^{\sigma \phi }\right) $. By arguing as
in the first case we get the same estimate (\ref{fordesigualfinal}) for $%
0=\delta _{n}=\lambda _{1}+\cdots +\lambda _{n}$, which is a contradiction.
\end{enumerate}

This proves Theorem \ref{teocont} for the group $\mathrm{Sl}\left( d,\mathbb{%
R}\right) $.

Now we consider a general semi-simple group $G$. As before let $\Sigma
=\{\alpha _{1},\ldots ,\alpha _{l}\}$ and $\Delta =\{\delta _{1},\ldots
,\delta _{l}\}$ be the simple system of roots and fundamental weights,
respectively. Given $\Theta \subset \Sigma $ let $\Delta _{\Theta }\subset
\Delta $ be the set of fundamental weights $\delta _{j}$ such that the root
with the same index $\alpha _{j}\in \Theta $. We put%
\[
\mathfrak{a}\left( \Theta \right) =\mathrm{span}\left( \Theta \right) \qquad
\mathfrak{a}_{\Theta }=\mathrm{span}\left( \Delta \setminus \Delta _{\Theta
}\right) .
\]%
These subspaces are orthogonal to each other, and since $\Theta \cup \left(
\Delta \setminus \Delta _{\Theta }\right) $ is a basis of $\mathfrak{a}%
^{\ast }$ \ we have $\mathfrak{a}=\mathfrak{a}\left( \Theta \right) \oplus
\mathfrak{a}_{\Theta }$.

The proof of continuity will be an easy consequence of the following
algebraic lemma.

\begin{lema}
If $\delta \in \Delta _{\Theta }$ then its coordinates with respect to $%
\Theta \cup \left( \Delta \setminus \Delta _{\Theta }\right) $ are
nonnegative.
\end{lema}

\begin{profe}
Let $\gamma _{1}$ and $\gamma _{2}$ be the orthogonal projections of $\delta
$ on $\mathfrak{a}\left( \Theta \right) $ and $\mathfrak{a}_{\Theta }$,
respectively. First we check that the coefficients of $\gamma _{1}$ with
respect to $\Theta $ are nonnegative. By definition there exists just one
root $\alpha \in \Theta $ such that $2\langle \alpha ,\delta \rangle
/\langle \alpha ,\alpha \rangle =1$ and $2\langle \beta ,\delta \rangle
/\langle \beta ,\beta \rangle =0$ if $\beta \neq \alpha $. But if $\beta \in
\Theta $ then $2\langle \beta ,\delta \rangle /\langle \beta ,\beta \rangle
=2\langle \beta ,\gamma _{1}\rangle /\langle \beta ,\beta \rangle $. Hence $%
\gamma _{1}$ \ is a fundamental weight for the root system defined by $%
\Theta $. Its coefficients with respect to $\Theta $ are the entries of the
inverse of the Cartan matrix, which are nonnegative.

As to $\gamma _{2}$ it is given by the mean
\[
\gamma _{2}=\frac{1}{|\mathcal{W}_{\Theta }|}\sum_{w\in \mathcal{W}_{\Theta
}}w\delta ,
\]%
because $\gamma _{2}$ is orthogonal to $\Theta $ and hence $w\gamma
_{2}=\gamma _{2}$ for every $w\in \mathcal{W}_{\Theta }$. Moreover, $\langle
\Theta \rangle $ is a root system in $\mathfrak{a}\left( \Theta \right) $
with Weyl group $\mathcal{W}_{\Theta }$ which has no fixed points in $%
\mathfrak{a}\left( \Theta \right) $ besides $0$. Hence the mean applied to $%
\gamma _{1}$ is $0$ since it is a fixed point.

Now if $\mathfrak{a}^{+}=\{\beta \in \mathfrak{a}^{\ast }:\forall \alpha \in
\Sigma ,~\langle \alpha ,\beta \rangle >0\}$ is the Weyl chamber in $%
\mathfrak{a}^{\ast }$ then the fundamental weight $\delta \in \mathrm{cl}%
\mathfrak{a}^{+}$. Hence 

\[
\gamma _{2}\in \mathcal{W}_{\Theta }\left(
\mathrm{cl}\mathfrak{a}^{+}\right) =\bigcup\limits_{w\in
\mathcal{W}_{\Theta }}w\left( \mathrm{cl} \mathfrak{a}^{+}\right)
\]%
because $\mathcal{W}_{\Theta }\left( \mathrm{cl}\mathfrak{a}^{+}\right) $ is
a cone. On the other hand if $\beta \notin \Theta $ and $\gamma \in \mathcal{%
W}_{\Theta }\left( \mathrm{cl}\mathfrak{a}^{+}\right) $ then $\langle \beta
,\gamma \rangle \geq 0$ (see e.g. \cite{smsec}, Lemma 7.5). But
\[
\gamma _{2}=\sum_{\beta \notin \Theta }\frac{2\langle \gamma _{2},\beta
\rangle }{\langle \beta ,\beta \rangle }\delta _{\beta }
\]%
where $\delta _{\beta }$ is the fundamental weight corresponding to $\beta $%
. Hence the coefficients of $\gamma _{2}$ are nonnegative, concluding the
proof.
\end{profe}

\vspace{12pt}%

\noindent%
\textbf{Proof of Theorem \ref{teocont} for }$G$ \textbf{semi-simple:} Let $%
\Theta =\Theta _{\mathrm{Ly}}\left( \phi \right) =\Theta _{\mathrm{Mo}}$ and
take $\delta \in \Delta _{\Theta }$. Then $\sigma \mapsto \delta \left( H_{%
\mathrm{Ly}}^{\sigma \phi }\right) $ is upper semi-continuous. Write
\[
\delta =\sum_{\alpha \in \Theta }a_{\alpha }\alpha +\sum_{\delta \in \Delta
\setminus \Delta _{\Theta }}b_{\delta }\delta .
\]%
with $a_{\alpha }\geq 0$ by the lemma. We have
\[
\delta \left( H_{\mathrm{Ly}}^{\sigma \phi }\right) =\sum_{\alpha \in \Theta
}a_{\alpha }\alpha \left( H_{\mathrm{Ly}}^{\sigma \phi }\right)
+\sum_{\lambda \in \Delta \setminus \Delta _{\Theta }}b_{\lambda }\lambda
\left( H_{\mathrm{Ly}}^{\sigma \phi }\right)
\]%
where the last sum is continuous by the differentiability result of \cite%
{ferrsm}. Hence the first sum is upper semi-continuous as well. The
assumption $\Theta _{\mathrm{Ly}}\left( \phi \right) =\Theta _{\mathrm{Mo}}$
implies that the first sum is zero at $\sigma =\mathrm{id}$. Since it is
nonnegative because $\alpha \left( H_{\mathrm{Ly}}^{\sigma \phi }\right)
\geq 0$ and $a_{\alpha }\geq 0$ we conclude that $\delta $ is continuous at $%
\mathrm{id}$, proving Theorem \ref{teocont}.

It remains to consider the reductive groups, which amounts to check
continuity of the central component defined in Section 3.3 of \cite{alvsm}.
The continuity of this component holds without any further assumption. This
is because this central component is given by an integral
\[
\int \mathfrak{a}_{\sigma }^{+}\left( 1,x\right) \nu \left( dx\right)
\]%
on the base space whose integrand is the time $1$ of a cocycle $\mathfrak{a}%
_{\sigma }^{+}\left( n,x\right) $ that depends continuously of $\sigma \in
\mathcal{G}$ (see \cite{alvsm} for the details).

\end{document}